\documentclass[11pt,a4paper, twoside]{article}
\usepackage{amssymb,amsfonts,amsmath,amsthm}
\usepackage[latin1]{inputenc}
\usepackage[english]{babel}
\usepackage[left=2.7cm,right=2.7cm,bottom=2.67cm,top=2.67cm]{geometry}
\usepackage{mathtools}
\usepackage{dsfont}
\usepackage{mathrsfs}
\usepackage{natbib}
\usepackage{enumerate}
\usepackage{color}
\usepackage{soul}
\usepackage[implicit=false]{hyperref}
\RequirePackage[OT1]{fontenc}
\usepackage{graphicx}
\usepackage{multirow}
\usepackage{multicol}
\usepackage{subfigure}
\usepackage{placeins}
\usepackage{rotating}
\usepackage[flushleft]{threeparttable}
\usepackage{booktabs}
\usepackage[normalem]{ulem}

\graphicspath{{figure/}}

\theoremstyle{plain}

\newtheorem{lema}{Lemma}

\theoremstyle{definition}

\newcommand{\E}{\mathds{E}}
\newcommand{\F}{\mathscr{F}}
\newcommand{\R}{\mathds{R}}

\newcommand{\Z}{\mathds{Z}}

\newcommand{\bs}[1]{\boldsymbol{#1}}

\long\def\sfootnote[#1]#2{\begingroup%
\def\thefootnote{\fnsymbol{footnote}}\footnote[#1]{#2}\endgroup}

\def\bfootnote{\xdef\@thefnmark{}\@footnotetext}

\begin{document}
\pagestyle{myheadings} % allow the use of headings
\markboth{Unit-Lindley Autoregressive Moving Average }{G. Pumi, D.H. Matsuoka and T.S. Prass}

\thispagestyle{empty}
{\centering
\Large{\bf A GARMA Framework for Unit-Bounded Time Series Based on the Unit-Lindley Distribution with Application to Renewable Energy Data}\vspace{.5cm}\\
\normalsize{ {\bf Guilherme Pumi${}^{\mathrm{a,}}$\sfootnote[1]{Corresponding author. This Version: \today},\let\thefootnote\relax\footnote{\hskip-.3cm$\phantom{s}^\mathrm{a}$Mathematics and Statistics Institute and Programa de P\'os-Gradua\c c\~ao em Estat\'istica - Universidade Federal do Rio Grande do Sul.
%- 9500,  Bento Gon\c calves Avenue - 91509-900, Porto Alegre - RS - Brazil.
} Danilo Hiroshi Matsuoka${}^\mathrm{a}$ and Taiane Schaedler Prass${}^\mathrm{a}$ 
 \\
\let\thefootnote\relax\footnote{E-mails: guilherme.pumi@ufrgs.br (Pumi), danilomatsuoka@gmail.com (Matsuoka), taiane.prass@ufrgs.br (Prass).}
\let\thefootnote\relax\footnote{ORCIDs: 0000-0002-6256-3170 (Pumi); 0000-0002-9744-8260 (Matsuoka); 0000-0003-3136-909X (Prass).}}\\
\vskip.3cm
}}

\begin{abstract}
The Unit-Lindley is a one-parameter family of distributions in $(0,1)$ obtained from an appropriate transformation of the Lindley distribution. In this work, we introduce a class of dynamical time series models for continuous random variables taking values in $(0,1)$ based on the Unit-Lindley distribution. The models pertaining to the proposed class are observation-driven ones for which, conditionally on a set of covariates, the random component is modeled by a Unit-Lindley distribution. The systematic component aims at modeling the conditional mean through a dynamical structure resembling the classical ARMA models. Parameter estimation in conducted using partial maximum likelihood, for which an asymptotic theory is available. Based on asymptotic results, the construction of confidence intervals, hypotheses testing, model selection, and forecasting can be carried on. A Monte Carlo simulation study is conducted to assess the finite sample performance of the proposed partial maximum likelihood approach.  Finally, an application considering forecasting of the proportion of net electricity generated by conventional hydroelectric power in the United States is presented. The application show the versatility of the proposed method compared to other benchmarks models in the literature.
\vspace{.2cm}

\vspace{.2cm}

\end{abstract}
\textbf{Statements and Declarations:} The authors declare that they have NO affiliations with or involvement in any organization or entity with any financial interests in the subject matter or materials discussed in this manuscript.
\section{Introduction}
% Applications:
% UWARMA: Manufacturing Capacity Utilization
% MARMA: useful water volume
% KARMA: relative humidity
% BARC: proportion of stocked hydroelectric energy
% Scher: proportion of stocked hydroelectric energy
% BARMA:  hidden unemployment due to substandard work conditions
% BARFIMA: RH
% BURARMA: proportion of stored energy
%
In environmental sciences, many key indicators and climate-related proportions evolve over time and are naturally bounded within the interval $(0,1)$, making models tailored for $(0,1)$-bounded data essential for accurate modeling, forecasting, and informed decision-making in the context of environmental variability and climate change. For instance, in this paper we consider the problem of modeling and forecasting the proportion of net electricity generated by conventional hydroelectric power in the United States. Hydroelectric power remains one of the most significant sources of renewable energy worldwide and plays a crucial role in the U.S. energy mix, consistently generating about three times more electricity than nuclear power. According to the time series data available from the U.S. Energy Information Administration (EIA), it was only in 2022 that the combined energy output of all other renewable sources (excluding nuclear) matched that of conventional hydro. Given this context, accurately modeling and forecasting the proportion of hydroelectric power in total electricity generation is of utmost importance for informed decision-making and long-term energy planning. Moreover, since this proportion is inherently affected by climatic and hydrological variability, producing reliable forecasts helps mitigate uncertainties and supports the optimization of energy distribution. Understanding how this share fluctuates over time is essential for managing energy resources effectively in an evolving and increasingly renewable-focused energy landscape.

The constrained nature of such time series hinders the applicability of traditional Gaussian-based models, calling for specialized statistical tools to approach the problem. To address this, the literature on non-Gaussian time series models defined on $(0,1)$ has grown significantly over the last decade. One of the most promising and studied approaches for modeling $(0,1)$-valued time series is the class of Generalized Autoregressive Moving Average (GARMA) models, which gained visibility following the seminal work of \cite{garma}, where the terminology was formalized and existing approaches were systematized under a unified framework. 

The broad applicability of these models promoted further methodological developments and extensions across various fields. In environmental sciences, for instance, \cite{barc} and \cite{burarma} proposed novel models for the proportion of stored energy, \cite{karma, barfima} focused on relative air humidity, and more recently, \cite{marma} introduced a GARMA model tailored to model useful water volume in reservoirs. Similarly, in econometrics, these models have been employed to analyze phenomena such as hidden unemployment due to substandard working conditions \citep{barma}, manufacturing capacity utilization \cite{uwarma}, the proportion of income allocated to consumption, industry market shares, and the fraction of weekly working hours, among others.

While the GARMA models introduced in \cite{garma} only considered random components belonging to the canonical exponential family, the literature has transcended this limitation in order to  accommodate  random components that do not belong to the exponential family, leading to the nowadays called GARMA-like models. This is the case of the  $\beta$ARMA, \citep{barma}, KARMA \citep{karma}, $\beta$ARFIMA \citep{barfima}, $\beta$ARC \citep{barc}, UBXII-ARMA \citep{burarma}, to cite a few and only considering models restricted to $(0,1)$. There are methodological advantages of considering a random component from the canonical exponential family. For instance, obtaining formulas for the score and the information matrix are considerably simplified, inference using partial likelihood is very well understood, and its theoretical properties, including asymptotic theory, are well developed \citep{Fokianos2004}. %But this is not a limitation as  many of the most commonly applied models for time series in $(0,1)$ are GARMA-like. 

GARMA models are a particular case of the broader family of observation-driven models of \cite{cox1981}, which are defined by two components. The first one, called the random component, is responsible for the conditional probability structure of the model, given the observed history. The second one is the systematic component, which is responsible for the model's dependence structure. In GARMA models, this is accomplished in a generalized linear model fashion: a linear structure is prescribed to a certain quantity of interest, such as the conditional mean, median, $\rho$-quantile, etc., through a link. 

In this work we propose the class of ULARMA (Unit-Lindley Autoregressive Moving Average) models - a GARMA model based on the Unit-Lindley distribution parameterized by its mean, as introduced in \cite{mazucheli}. The Unit-Lindley is a versatile uniparametric distribution supported in $(0,1)$ belonging to the canonical exponential family. The systematic component is prescribed through an ARMA-like set of difference equations (see \eqref{eq2}), which can also contain a set of exogenous covariates. This structure is considered, for instance, in the $\beta$ARMA of \cite{barma}, based on the Beta distribution parametrized by its mean; the KARMA of \cite{karma}, which is based on the standard Kumaraswamy distribution parametrized by the median; and the UWARMA of \cite{uwarma} and UBXII-ARMA \citep{burarma}, based on the Unit-Weibull and Unit Burr XII distributions, respectively, both parametrized by their $\rho$th quantile. The systematic component may assume other forms, such as to include long-range dependence as in \cite{barfima} and \cite{symarfima}, generalized autoregressive score specifications as in \cite{rgas} or even using chaotic specifications as in \cite{barc}.

Inference for the ULARMA model is carried out using Partial Maximum Likelihood Estimation (PMLE), for which we derive closed-form expressions for both the partial score vector and the conditional information matrix. Building on this, we develop a comprehensive framework for large-sample inference, including hypothesis testing and forecasting. In particular, we introduce a bootstrap-based method to construct out-of-sample prediction intervals. To evaluate model adequacy, we discuss diagnostic tools such as residual analysis and goodness-of-fit procedures tailored to the proposed setting. The finite-sample behavior of the PMLE method is investigated through a Monte Carlo simulation study, focusing on point estimation accuracy, the performance of goodness-of-fit tests, and the joint behavior of the estimated parameters.

Given that the Unit-Lindley distribution is uniparametric and follows the systematic structure of several other more parameterized models, one fair question is why it is important to introduce and study a potentially simpler model. We identify three main reasons. First, the Unit-Lindley distribution belongs to the canonical exponential family, inheriting its desirable properties and benefiting from a well-established general theory, which simplifies foundational calculations. Second, as consequence, when inference is based on partial maximum likelihood for the ULARMA model, we obtain a solid asymptotic framework that supports essential inferential tools, such as the construction of (asymptotic) test statistics and confidence intervals. Third -- and most importantly -- despite its simplicity, the ULARMA model is highly versatile, as demonstrated in our empirical study. As mentioned before, we address the problem of modeling and forecasting the proportion of net electricity generated by conventional hydroelectric power in the United States using time series data from 2001 to 2025. We compare the forecasting performance of the ULARMA model with two benchmark models, namely the KARMA and $\beta$ARMA. The results are noteworthy: the ULARMA model fits the data more parsimoniously and yields more accurate 12-month-ahead forecasts than its competitors. This exemplifies the potential of a ``simpler'' model when applied in the appropriate context.

\section{The ULARMA model}
The Unit-Lindley distribution was introduced by \cite{mazucheli} upon considering a transformation of the standard Lindley distribution \citep{lindley} to the unit. Its probability density function is defined by
\begin{equation}\label{eq1}
%    f(y;\mu)=\frac{(1-\mu)^2}{\mu(1-y)^3}\exp\bigg( \frac{\mu-1}{\mu}\frac{y}{1-y}\bigg), \ 0<y<1,
    f(y;\mu)=\frac{(1-\mu)^2}{\mu(1-y)^3}\exp\bigg\{ \frac{y(\mu-1)}{\mu(1-y)}\bigg\} I(0<y<1),
\end{equation}
where $\mu\in(0,1)$. We use the notation $X \sim \mathrm{UL}(\mu)$ to say that a random variable $X$ has the Unit-Lindley distribution with parameter $\mu$. 

It easy to see that, under parametrization \eqref{eq1}, $\E(X)=\mu$. This relation allowed \cite{mazucheli} to proposed a GLM model based on the Unit-Lindley distribution. In this work we extend their approach by proposing a GARMA model based on the Unit-Lindley distribution. Let $\{Y_t\}_{t\in\Z}$ be a stochastic process taking values in $(0,1)$ and let $\{\bs X_t\}_{t\in\Z}$ be a set of $r$-dimensional exogenous covariates to be included in the model. These can be either random or deterministic and time-dependent, or any combination of these.  Without loss of generality, throughout the text it is assumed that $X_{tk}$ denotes the value of the $k$th covariate at time $t$, for deterministic covariates, and at time $t -1$, for random ones,  for $1 \leq k \leq r$. Let $\F_{t}$ denote the information (sigma-field) available to the observer at time $t$, that is  $\F_{t}:=\sigma\{\bs X_{t+1}, Y_{t},\bs X_{t},Y_{t-1},\cdots\}$ The proposed Unit-Lindley Autoregressive Moving Average (ULARMA) is a class of observation-driven models for which the random component is implicitly defined by assigning  $Y_t|\F_{t-1}\sim \mathrm{UL}(\mu_t)$. Upon observing that  $\mu_t:=\E(Y_t|\F_{t-1})$, we follow the GARMA approach by setting
\begin{equation}\label{eq2}
\eta_t:=g(\mu_t)= \alpha + \bs X_t^\prime \bs\beta + \sum_{i=1}^p \phi_i \big[ g(Y_{t-i})-\bs X_{t-i}^\prime \bs \beta\big] + \sum_{j=1}^q \theta_j r_{t-j},
\end{equation}
where $\eta_t$ is the linear predictor, $g:(0,1)\rightarrow \R$ is a twice differentiable link, $\alpha$ is an intercept, $\bs\beta=(\beta_1, \cdots,\beta_r)^\prime$ is the parameter vector related to the covariates, $\bs\phi=(\phi_1,\cdots,\phi_p)^\prime$ and $\bs\theta=(\theta_1,\cdots,\theta_q)^\prime$ are the AR and MA coefficients, respectively. The error term  in \eqref{eq2} is defined in a recursive fashion by setting $r_t:=g(Y_t)-g(\mu_t)$. The proposed class of models, hereafter denoted as ULARMA$(p,q)$, is defined by setting $Y_t|\F_{t-1}\sim \mathrm{UL}(\mu_t)$, with $\mu_t$ specified by \eqref{eq2}.

\section{Parameter estimation}

In this section we propose the use of the partial maximum likelihood approach to parameter estimation in the context of ULARMA models. Let $\{Y_t\}_{t\in\Z}$ be an ULARMA$(p,q)$ model with associated $r$-dimensional covariates $\{\bs X_t\}_{t\in\Z}$. Let $\bs\gamma:=(\alpha,\bs\beta,\bs\phi,\bs\theta)^\prime \in \Omega$, where $\Omega\subset\R^{r+p+q+1}$ denotes the parameter space. Given a sample $\{(Y_t,\bs X_t)\}_{t=1}^n$, the partial log-likelihood function is given by 
\begin{equation*}%\label{logL}
\ell(\bs \gamma)= \sum_{t=1}^{n} \ell_t(\bs \gamma),
\end{equation*}
where
\begin{align*}
%\ell_t(\bs \gamma)&:= 2\ln(1-\mu_t)-\ln(\mu_t)-3\ln(1-Y_t)+\frac{\mu_t-1}{\mu_t}\frac{Y_t}{1-Y_t}.
\ell_t(\bs \gamma):= 2\ln(1-\mu_t)-\ln(\mu_t)-3\ln(1-Y_t)+\frac{Y_t(\mu_t-1)}{\mu_t(1-Y_t)}.
\end{align*}
 The partial maximum likelihood estimator (PMLE) of $\bs\gamma$ is given by
\begin{equation*}
\hat{\bs\gamma}=\underset{\bs\gamma\in\Omega}{\mathrm{argmax}}\big\{\ell(\bs\gamma)\big\}.
\end{equation*}
\subsection{The partial score vector}
 The partial score vector is defined as $\frac{\partial \ell_t(\bs\gamma)}{\partial \bs\gamma}$. For $j\in\{1,\dotsc,p+q+r+1\}$, we have that
\begin{align*}%\label{dldnu}
\frac{\partial \ell_t(\bs \gamma)}{\partial \gamma_j} 
&=\frac{\partial \ell_t(\bs \gamma)}{\partial \mu_t}\frac{\partial \mu_t}{\partial \eta_t}\frac{\partial \eta_t}{\partial \gamma_j}=\frac1{g'(\mu_t)}\bigg[-\frac{2}{1-\mu_t}-\frac{1}{\mu_t}+\frac{Y_t}{\mu_t^2(1-Y_t)}\bigg] \frac{\partial \eta_t}{\partial \gamma_j}
\end{align*}
and 
\begin{align}\label{equas}
\frac{\partial \eta_t}{\partial \alpha}&=1 - \sum_{j=1}^q \theta_j \frac{\partial \eta_{t-j}}{\partial \alpha},\qquad\qquad
\frac{\partial \eta_t}{\partial \beta_l}= X_{tl}-\sum_{i=1}^{p} \phi_iX_{(t-i)l}- \sum_{j=1}^q \theta_j \frac{\partial \eta_{t-j}}{\partial \beta_l},\nonumber\\
\frac{\partial \eta_t}{\partial \phi_k}&= g({Y}_{t-1})-\bs X_{t-1}'\bs\beta - \sum_{j=1}^q \theta_j \frac{\partial \eta_{t-j}}{\partial \phi_k},\qquad\mbox{and}\qquad
\frac{\partial \eta_t}{\partial \theta_s}=r_{t-j} - \sum_{i=1}^q \theta_i \frac{\partial \eta_{t-i}}{\partial \theta_s},
\end{align}
for $l\in\{1,\cdots,r\}$, $k\in\{1,\cdots,p\}$ and $s\in\{1,\cdots,q\}$, where $X_{tl}$ denotes the $l$-th component of $\bs X_t$. 
The score vector $U(\bs\gamma)$ can be conveniently written in matrix form as
\begin{equation*}%\label{eq:score}
U(\bs\gamma) = D' T\bs h  
\end{equation*}
where $D$ is the $n\times(p+q+r+1)$ matrix whose $(i,j)$th element is given by $[D]_{i,j} = \partial \eta_{i}/\partial \gamma_j$, $T$ is a diagonal matrix and $\bs h$ a vector given respectively by
\begin{align*}
T:=\mathrm{diag}\biggl\{\frac{1}{g_1'(\mu_1)},\cdots,\frac{1}{g_1'(\mu_n)}\biggr\} \qquad \mbox{and} \qquad \bs{h} := \bigg(\frac{\partial \ell_1(\bs\gamma)}{\partial \mu_1}, \cdots, \frac{\partial \ell_n(\bs\gamma)}{\partial \mu_n}\bigg)'.
\end{align*}
The PMLE of $\bs\gamma$ if it exists, it is obtained as a solution of the so-called normal equations, given by the system $U(\gamma)=\mathbf{0}$, where $\mathbf{0}$ is the null vector in $\R^{p+q+r+1}$. However, the normal equations cannot be solved analytically. In this case, we have to resort to numerical optimization to approximate the PMLE.

\subsection{Conditional information matrix}

In this section we derive the Fisher conditional information matrix, which will be useful later on deriving the asymptotic properties of the partial maximum likelihood estimator for the proposed model. Let 
\begin{equation*}
H_t(\bs \gamma) := -\frac{\partial^2\ell_t(\bs\gamma)}{\partial \bs \gamma \partial \bs \gamma'},\quad\mbox{and}\quad
H(\bs \gamma) := -\frac{\partial^2\ell(\bs\gamma)}{\partial \bs \gamma \partial \bs \gamma'}  =  -\sum_{t=1}^{n}\frac{\partial^2\ell_t(\bs\gamma)}{\partial \bs \gamma \partial \bs \gamma'} = \sum_{t=1}^nH_t(\bs \gamma).
\end{equation*}
Notice that $H(\bs \gamma)$ and $\ell(\bs\gamma)$ both depend on $n$. However, for simplicity and since no confusion will arise, we shall drop the dependence on $n$ on the notation. Let $I_n(\bs \gamma) := \E(H(\bs \gamma))$ be the information matrix corresponding to the sample of size $n$ and let
\begin{equation*}
   I^{(n)}(\bs \gamma) := -\frac{1}{n}\sum_{t = 1}^n \E \bigg(\frac{\partial^2\ell_t(\bs \gamma)}{\partial \bs\gamma \partial \bs \gamma'} \bigg)= -\frac{1}{n}\E \bigg(\frac{\partial^2\ell(\bs \gamma)}{\partial \bs \gamma \partial \bs \gamma'} \bigg),
\end{equation*}
so that $I_n(\bs \gamma ) = nI^{(n)}(\bs \gamma)$. Now, observe that
\begin{equation*}
 I^{(n)}(\bs \gamma) = -\frac{1}{n}\sum_{t = 1}^n \E \bigg( \E\bigg(\frac{\partial^2\ell_t(\bs \gamma)}{\partial \bs\gamma \partial \bs \gamma'} \bigg| \F_{t-1}\bigg)\bigg)  = \frac{1}{n}\E\big(K_n(\bs \gamma)\big),
\end{equation*}
with
\begin{equation*}
 K_n(\bs \gamma) : = -\sum_{t = 1}^n \E \bigg(\frac{\partial^2\ell_t(\bs \gamma)}{\partial \bs\gamma \partial \bs\gamma'} \Big| \F_{t-1}\bigg).
\end{equation*}
The matrix $K_n(\bs \gamma)$ is known as the conditional information matrix corresponding to the sample of size $n$. 
Under some regularity conditions presented later,
\begin{equation}\label{ah}
 \frac{1}{n}H(\bs \gamma) - I^{(n)}(\bs \gamma) \overset{P}{\longrightarrow} 0 \quad \mbox{and} \quad    \frac{1}{n}K_n(\bs \gamma) - I^{(n)}(\bs \gamma) \overset{P}{\longrightarrow} 0, \quad \mbox{as} \quad n\to \infty.
\end{equation}
Furthermore, $I^{(n)}(\bs \gamma) \, {\longrightarrow} \,  I(\bs \gamma)$, where
\begin{equation}\label{I}
I(\bs \gamma) = \lim_{n\to \infty} I^{(n)}(\bs \gamma) = \lim_{n\to \infty} -\frac{1}{n}\E \bigg(\frac{\partial^2\ell(\bs \gamma)}{\partial \bs \gamma \partial \bs \gamma'} \bigg),
\end{equation}
which is the analogous of the $I_1(\bs \gamma)$ matrix for i.i.d. samples. We shall derive $K_n(\bs\gamma)$ in closed form. First notice that
\begin{align*}
\frac{\partial^2\ell_t(\bs \gamma)}{\partial \gamma_i \partial \gamma_j} &= \frac{\partial}{\partial \mu_t}
\left( \frac{\partial \ell_t(\bs \gamma)}{\partial \mu_t}\frac{\partial \mu_t}{\partial \eta_t}\frac{\partial \eta_t}{\partial \gamma_j}\right)
\frac{\partial \mu_t}{\partial \eta_t} \frac{\partial \eta_t}{\partial \gamma_i} \\
&=  \left[ \frac{\partial^2 \ell_t(\bs \gamma)}{\partial \mu_t^2}\frac{\partial \mu_t}{\partial \eta_t}\frac{\partial \eta_t}{\partial \gamma_j}
+ \frac{\partial \ell_t(\bs \gamma)}{\partial \mu_t}\frac{\partial}{\partial \mu_t}\left(\frac{\partial \mu_t}{\partial \eta_t}\frac{\partial \eta_t}{\partial \gamma_j} \right) \right]\frac{\partial \mu_t}{\partial \eta_t} \frac{\partial \eta_t}{\partial \gamma_i} \,.
\end{align*} 
Observe that, by Lemma \ref{lema} in Appendix A, $\E\Big(\frac{\partial \ell_t(\bs\gamma)}{\partial \mu_t} \big| \F _{t-1}\Big)=0$.  The $\F_{t-1}$-measurability of $\mu_t$ and $\eta_t$, implies that
\begin{equation*}
\E\bigg(\frac{\partial^2\ell_t(\bs\gamma)}{\partial \gamma_i \partial \gamma_j}\Big|\F_{t-1}\bigg)=  \E\bigg(\frac{\partial^2 \ell_t(\bs\gamma)}{\partial \mu_t^2}\Big|\F_{t-1}\bigg)\bigg[\frac{\partial \mu_t}{\partial \mu_t}\frac{\partial \mu_t}{\partial \eta_t}\bigg]^2 \frac{\partial \eta_t}{\partial \gamma_i}\frac{\partial \eta_t}{\partial \gamma_j},
\end{equation*}
where $\frac{\partial \eta_t}{\partial \gamma_k}$ is given in \eqref{equas}. Hence, we only need to obtain $\E\Big(\frac{\partial^2 \ell_t(\bs\gamma)}{\partial \mu_t^2}\big|\F_{t-1}\Big)$. We have that
\begin{align*}
\frac{\partial^2 \ell_t(\bs \gamma)}{\partial \mu_t^2}&=\frac{\partial}{\partial \mu_t}\bigg[-\frac{2}{1-\mu_t}-\frac{1}{\mu_t}+\frac{Y_t}{\mu_t^2(1-Y_t)}\bigg]=-\frac{2}{(1-\mu_t)^2}+\frac{1}{\mu_t^2}-\frac{2Y_t}{\mu_t^3(1-Y_t)}.
\end{align*}
Now, after some algebra from the $\F_{t-1}$-measurability of $\mu_t$ and Lemma \ref{lema}, we obtain
\begin{align*}
\E\bigg(\frac{\partial^2 \ell_t(\bs \gamma)}{\partial \mu_t^2}\Big|\F_{t-1}\bigg)&=-\frac{2}{(1-\mu_t)^2}+\frac{1}{\mu_t^2}-\frac{2}{\mu_t^3}\E\bigg(\frac{Y_t}{1-Y_t}\bigg |\mathcal{F}_{t-1}\bigg)%\\
%&=-\frac{2}{(1-\mu_t)^2}+\frac{1}{\mu_t^2}-\frac{2}{\mu_t^3}\frac{\mu_t^2+\mu_t}{1-\mu_t}
%\\ &=\frac{-2\mu_t+(1-\mu_t)^2-2(1+\mu_t)(1-\mu_t)}{\mu_t^2(1-\mu_t)^2}\\&
=\frac{(1-\mu_t)^2-2}{\mu_t^2(1-\mu_t)^2}.
\end{align*}
By letting $E_{\mu}$ be the $n\times n$ diagonal matrix for which the $k$th diagonal elements is given by
\begin{equation*} 
[E_\mu]_{k,k}:=-\E\bigg(\frac{\partial^2 \ell_k(\bs\gamma)}{\partial \mu_k^2}\Big|\F_{k-1}\bigg)=\frac{2-(1-\mu_k)^2}{\mu_k^2(1-\mu_k)^2},
\end{equation*}
and $D_{\bs\gamma}$ and $T$ as before, we obtain
\begin{equation}\label{kn}
K_n(\bs\gamma)=D_{\bs\gamma}'TE_{\mu}TD_{\bs\gamma}.
\end{equation}
%%%%%%%%%%%%%%%%%%%%%%%%%%%%%%%%%%%%%%%#
%%%%%%%%%%%%%%%%%%%%%%%%%%%%%%%%%%%%%%%%%%%%%%%%%%
%%%%%%%%%%%%%%%%%%%%%%%%%%%%%%%%%%%%%%%%%%%%%%%%%%

\section{Large sample inference} \label{lsi}

The asymptotic properties of the PMLE for GARMA-type models, when the response distribution belongs to the canonical exponential family, have been rigorously developed in \cite{Fokianos2004}. These results apply directly to the ULARMA model proposed in this work.

For notational clarity and without loss of generality, consider specification \eqref{eq2} without covariates. Let $Y_1, \cdots, Y_n$ be a realization from a ULARMA$(p,q)$ process, and define the $(p+q+1)$-dimensional covariate vector as
\[ \bs Z_{t-1} := \big(1, g(Y_{t-1}), \cdots, g(Y_{t-p}), r_{t-1}, \cdots, r_{t-q} \big)', \]
so that the systematic component can be expressed compactly as $\eta_t = \bs Z_{t-1}'\bs\gamma$, with $\bs\gamma = (\alpha, \bs\phi', \bs\theta')'$. The regularity conditions ensuring the consistency and asymptotic normality of the PMLE $\hat{\bs\gamma}$ in this setting are as follows:
\begin{enumerate}
    \item The true parameter $\bs\gamma_0$ lies in an open set $\Omega \subseteq \R^{p+q+1}$, and the covariates $\bs Z_{t-1}$ are almost surely contained within a compact subset $\Gamma \subset \R^{p+q+1}$, such that
    \[P \left( \sum_{t=1}^n \bs Z_{t-1} \bs Z_{t-1}' > 0 \right) = 1. \]
    \item The link function $g$ is twice continuously differentiable, with inverse $g^{-1}$ satisfying $\partial g^{-1}(x)/\partial x \neq 0$. Furthermore, $\bs Z_{t-1}'\bs\gamma$ lies almost surely in the domain of $g^{-1}$ for all $\bs\gamma \in \Omega$ and all $t$.
    \item There exists a probability measure $\lambda$ on $\R^{p+q+1}$ such that the matrix  $\displaystyle{\int_{\R^{p+q+1}}} \bs v\bs v' \lambda(d\bs v)$ is positive definite, and the empirical measure satisfies
    \[\frac{1}{n} \sum_{t=1}^n \mathbb{I}(\bs Z_{t-1} \in A) \overset{P}{\longrightarrow} \lambda(A), \qquad \text{as } n \to \infty,  \]
    for every measurable set $A$.
\end{enumerate}
These conditions and their implications are discussed in detail in Section 5 of \cite{Fokianos2004}. In particular, condition 3 ensures a form of ergodicity: for any continuous and bounded function $h$ defined on $\Gamma$,
\[\frac{1}{n} \sum_{t=1}^n h(\bs Z_{t-1}) \overset{P}{\longrightarrow} \int_{\R^{p+q+1}} h(\bs v) \, \lambda(d\bs v),\]
which implies \eqref{ah}, with $I(\bs\gamma_0)$ positive definite and, thus, invertible. Under these regularity conditions, an almost surely unique PMLE $\hat{\bs\gamma}$ exists and satisfies $\hat{\bs\gamma} \overset{P}{\longrightarrow} \bs\gamma_0$, as $n \to \infty$.

Moreover, the estimator is asymptotically normal
\begin{equation} \label{clt}
\sqrt{n} (\hat{\bs\gamma} - \bs\gamma_0) \overset{d}{\longrightarrow} \mathcal{N}_{p+q+1} \big( \bs 0, I^{-1}(\bs\gamma_0) \big),
\end{equation}
where $I(\bs\gamma_0)$ is defined in \eqref{I}, and $\mathcal{N}_m(\bs 0, \Sigma)$ denotes the $m$-variate normal distribution with mean vector $\bs 0 \in \R^m$ and covariance matrix $\Sigma$. The proof of these results relies on analyzing the asymptotic behavior of the conditional information matrix and applying a central limit theorem to the properly normalized partial score vector $U(\bs\gamma)$. For a detailed derivation, see \cite{Fokianos1998, Fokianos2004}.

\subsection{Hypothesis Testing} \label{s:hp}

Let $Y_1, \cdots, Y_n$ be a sample from a ULARMA$(p,q)$ model that includes a set of $r$-dimensional covariates $\bs X_1, \cdots, \bs X_n$. Denote the true parameter vector by $\bs\gamma_0 := (\gamma_1^0, \cdots, \gamma_{p+q+r+1}^0)'$, and let $\hat{\bs\gamma} := (\hat\gamma_1, \cdots, \hat\gamma_{p+q+r+1})'$ be the corresponding PMLE. Let $K_n(\hat{\bs\gamma})$ be the empirical counterpart of the conditional information matrix as defined in \eqref{kn}, evaluated at the estimate $\hat{\bs\gamma}$. Denote by $K_n(\hat{\bs\gamma})^{jj}$ the $j$th diagonal element of the inverse matrix $K_n(\hat{\bs\gamma})^{-1}$. Under the regularity conditions established in Section~\ref{lsi}, the asymptotic normality result in \eqref{clt} enables the construction of large-sample hypothesis tests. For example, to test the null hypothesis $H_0: \gamma_j = \gamma_j^\ast,$ for a given value $\gamma_j^\ast$, one may use the Wald-type statistic
\[z = \frac{\hat\gamma_j - \gamma_j^\ast}{\sqrt{K_n(\hat{\bs\gamma})^{jj}}},\]
which, under $H_0$, is approximately distributed as a standard normal variable for sufficiently large $n$.

This approach generalizes to joint tests involving subsets of parameters and can be used to assess the significance of autoregressive, moving average, or regression components. Other classical testing procedures -- such as likelihood ratio and score tests -- can also be adapted to this framework and retain their usual asymptotic distributions, as in the case of independent observations. For detailed derivations and further discussion, see \cite{Fahrmeir1987} and Section 6 of \cite{Fokianos2004}.

\subsection{Forecasting}\label{msf}
Let $Y_1,\cdots, Y_n$ be a sample of a ULARMA$(p,q)$ model with associated covariates $\bs X_1,\cdots,\bs X_n$. $h$-step ahead forecasts  $\hat y_{n+1}, \cdots, \hat y_{n+h}$ can be obtained recursively from an adapted sample version of \eqref{eq2}. Observe that, in the presence of covariates in the model, the values of $\bs X_{n+1},\cdots,\bs X_{n+h}$ are required to construct a sample version of \eqref{eq2}.  The simplest case is when the covariates are deterministic or, to a less extent, predetermined. However, if the covariates are random, we assume that $\bs X_{n+1},\cdots,\bs X_{n+h}$ are available or can be obtained (by forecasting, for instance). Let $\hat{\bs\gamma}$ denote the PMLE estimated from the sample, in view of \eqref{eq2}, starting at $t=1$, we recursively obtain 
\begin{equation}\label{rec}
\hat \eta_{t}  = \hat\alpha + \hat{\bs X}_t'\hat{\bs\beta} + \sum_{i=1}^p  \hat\phi_i\bigl[g_{2}(\hat Y_{t-i})-\hat{\bs X}_{t-i}'\hat{\bs\beta}\bigr] + \sum_{k=1}^q \hat \theta_k \hat r_{t-k},
\end{equation}
with $\hat \mu_t  = g_1^{-1}(\hat \eta_t)$, for $t \geq 1$, $\hat r_t = (g(\hat Y_t) - \hat \eta_t)I(1 \leq t \leq n)$, 
\begin{equation*}%\label{cases}
    \hat Y_t = \begin{cases}
    g^{-1}(0), & p > 0, \ t < 1,\\
   Y_t, & 1 \leq  t \leq n,\\
  \hat \mu_t, & t > n,
   \end{cases}\qquad\mbox{and}\qquad
   \hat{\bs X}_t = \begin{cases}
\frac{1}{p}\sum\limits_{i = 1}^p \bs X_i, & p >0, \ t < 1,\\
\bs X_t, & t \geq 1.
\end{cases}
\end{equation*}
Observe that the sequence $\hat\mu_{1},\cdots,\hat\mu_{n}$ corresponds to the in-sample forecasted values, whereas by setting $\hat y_{n+k}=\hat\mu_{n+k}$ for $k\in\{1,\cdots,h\}$, we obtain $h$-step ahead forecasts.
\subsection{Residuals and goodness-of-fit}\label{gof}
Goodness-of-fit tests and residual analysis in the context of ULARMA models is carried on differently from linear models such as the classical ARMA. The main difference is that the error term in specification \eqref{eq2} is constructed iteratively, so that no  information about its distribution is available and we need to resort to alternative methods to properly carry on residual analysis. For instance, Observing that since $\mu_t=\E(Y_t|\F_{t-1})$, the simple residuals defined by $e_t:=Y_t-\mu_t$ are such that $\{(e_t,\F_{t-1})\}_{t\in\Z}$ is a martingale difference sequence. Hence, if the ULARMA model is well specified, $\hat e_t:=Y_t-\hat\mu_t$ is expected to behave as a martingale difference with respect to $\F_{t-1}$. This can be tested using a martingale difference test. 
Another commonly used approach is based on the so-called quantile residuals, defined as
\begin{equation}\label{qres}
e_t^{(q)}:=\Phi^{-1}\big(F(Y_t|\F_{t-1})\big),
\end{equation}
%Quantile \[Q(p \mid \theta) = \frac{1 + \theta + W_{-1} \left( (1+\theta)(p-1) e^{-(1+\theta)} \right)}{1 + W_{-1} \left( (1+\theta)(p-1) e^{-(1+\theta)} \right)},\]
where $F(\cdot|\F_{t-1})$ denotes the cumulative distribution function associated with the model's random component and $\Phi$ denote the standard normal cumulative distribution function. In the present setting, the quantile residuals obtained from \eqref{qres} after plugin-in the PMLE on \eqref{eq2}, asymptotically follow a multivariate standard normal distribution when the model is correctly specified \citep[see lemma 2.1 in][]{kall}. %These goodness-of-fit procedures will be further explored in the simulations.
\subsection{Confidence Intervals and Model Selection}\label{spi}
Under the assumptions presented in Section \ref{lsi} the approximation $\sqrt{K_n(\hat{\bs\gamma})^{jj}}(\hat\gamma_j-\gamma_j^0)\approx N(0,1)$ holds, for large enough $n$, in view of \eqref{clt}. Level $\delta$ approximate confidence intervals for $\gamma_j^0$ can be obtained straightforwardly as $\hat\gamma_j\pm z_{1-\delta/2}/\sqrt{K_n(\hat{\bs\gamma})^{jj}}$, where $z_{1-\delta/2}$ is the $(1-\delta/2)$th quantile from a standard normal distribution.  In view of \eqref{ah} and \eqref{clt}, the delta method allows for the construction of approximate level $\delta$ in-sample forecasting interval through \citep{Fokianos2004}
\begin{equation*}
CI(\mu_t;\delta) = \hat\mu_t \pm z_{1-\frac\delta2}\sqrt{\frac{\bs Z_{t-1}'K_n(\hat{\bs\gamma})^{-1}\bs Z_{t-1}}{ng'(\hat\mu_t)^2}}.
\end{equation*}
Out-of-sample forecasting intervals can be obtained using an approach closely related to parametric bootstrap. Let $y_1,\cdots,y_n$ be a sample from an ULARMA$(p,q)$ model with estimated PMLE $\hat{\bs\gamma}$. To obtain $h$-step ahead forecasting intervals for $Y_t$, the proposed bootstrap approach relies on generating $B$ bootstrap samples for $Y_{n+1}, \cdots,Y_{n+h}$ from an ULARMA$(p,q)$ model considering $\hat{\bs\gamma}$ as parameter and updating the values of $\mu_t$ as the procedure evolves. 

To be more specific, we start by reconstructing the sequences $\hat{\mu}_1,\cdots,\hat{\mu}_n$ and $\hat r_{1}, \cdots, \hat r_n$ using the PMLE $\hat{\bs\gamma}$ through \eqref{rec}, as outlined in Section \ref{msf}. For each $m\in\{1,\cdots,B\}$, we start by setting $k=1$, and sampling $\hat y_{n+1}^{(m)}$ from a Unit-Lindley distribution with parameter $\hat\mu_{n+1}$ which can be obtained from \eqref{rec} using the available information. Next, we obtain $\hat r_{n+1}^{(m)} = g(\hat y_{n+1}^{(m)})-g(\hat \mu^{(m)}_{n+1})$ and perform the following two steps sequentially, for $k \in\{2,\cdots,h\}$:
\begin{enumerate}
\item Update $\hat\mu_{n+k}^{(m)}$ through \eqref{rec} using the augmented sample $y_1,\cdots, y_n,\hat y_{n+1}^{(m)}, \cdots\hat y_{n+k-1}^{(m)}$ and error terms $\hat r_1,\cdots, \hat r_n,$ $\hat r_{n+1}^{(m)}, \cdots,\hat r_{n+k-1}^{(m)}$.%, with $\hat y_{n}^{(m)} := y_n$ and $\hat r_{n}^{(m)}:= \hat r_{n}$.
\item  Sample $\hat y^{(m)}_{n+k}$ from a Unit-Lindley distribution with parameter $\hat\mu_{n+k}^{(m)}$ and update $\hat r_{n+k}^{(m)} = g(\hat y^{(m)}_{n+k})-g(\hat \mu^{(m)}_{n+k})$.
\end{enumerate}
After iterating these steps, we end up with a collection $\big\{\hat y^{(m)}_{n+1},\cdots, \hat y^{(m)}_{n+h}\big\}_{m=1}^B$ of $B$ bootstrapped samples. These can be used to approximate the distribution of $(Y_{n+1}, \cdots,Y_{n+h})$. For each $k\in\{1,\cdots,h\}$, an approximate confidence interval for $Y_{n+k}$ is obtained from the $(1-\delta/2)$th and $\delta/2$th sample quantiles calculated from $\hat y^{(1)}_{n+k}, \cdots,\hat y^{(B)}_{n+k}$. 

Finally, model selection in the context of GARMA model may be conducted using information criteria such as the AIC, BIC and HQC, which are calculated as per usual based on the maximized partial likelihood. Bayesian approaches via Reversible Jump Markov Chains have also been considered in the literature \citep{casarin, katerine}, but we shall not delve into this matter in this paper.
%%%%%%%%%%%%%%%%%%%%%%%%%%%%%%%%%%%%%%%%%%%%%%%%%%
%%%%%%%%%%%%%%%%%%%%%%%%%%%%%%%%%%%%%%%%%%%%%%%%%%
%%%%%%%%%%%%%%%%%%%%%%%%%%%%%%%%%%%%%%%%%%%%%%%%%%
\section{Monte Carlo Simulation}
In this section, we investigate the finite sample performance of the proposed PMLE method for parameter estimation in ULARMA models. Our objective is to examine  point estimation and residual analysis. The latter will be assessed using the methods outlined in Section \ref{gof}. The simulations were conducted using R version 4.3.1 \citep{R}.

Firstly, we consider the simple residuals, denoted as $\hat{e}_t = y_t - \hat{\mu}_t$, where $\hat{\mu}_t$ is derived from the PMLE. Under the assumption that the model is correctly specified, $\hat{e}_t$ should approximately follow a martingale difference process relative to the history of the process. This can be evaluated using various martingale difference tests, based on multiple approaches. In this study, we employ the wild bootstrap automatic variance ratio test, as proposed by \cite{kim2009}, and the Dom\'inguez-Lobato test \citep{DL}. The finite sample performance of this and other methods is further discussed in \cite{charles}. The second approach involves the use of quantile residuals, which, assuming correct model specification, should exhibit behavior consistent with a standard normal distribution.
\subsection{Point estimation exercise}
\begin{figure}[ht!]
\centering
\mbox{\includegraphics[width=\textwidth]{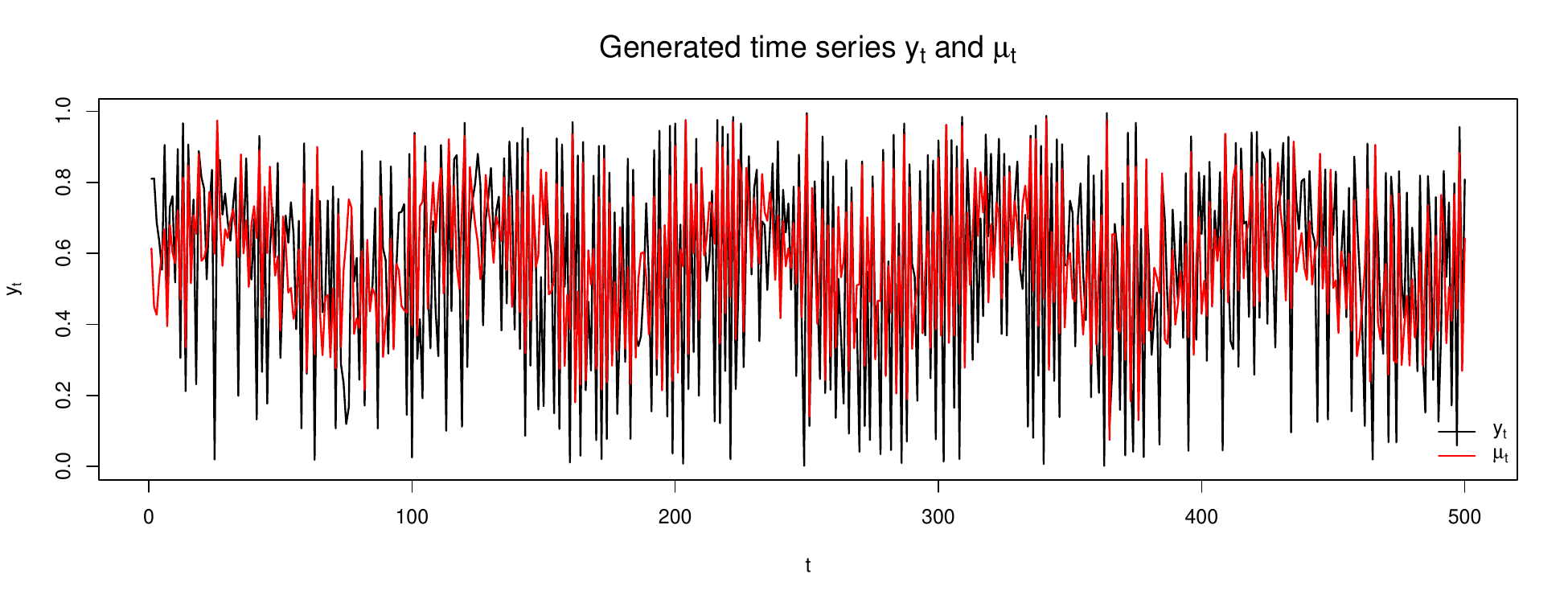}}
\caption{A typical example of a time series considered in the simulation study. The plot was generated considering $n=500$, $\alpha = 0.5$, $\beta=0.5$, $\phi=-0.4$, $\theta=-0.2$.} \label{ref}
\end{figure}

We start by considering the finite sample performance of the PMLE in the context of ULARMA models. We generate samples of size $n\in\{100,200,500\}$ from an ULARMA$(1,1)$ model with parameters $\alpha\in\{0.5,1\}$, $(\phi,\theta) \in\big\{(0.2,-0.8),(-0.8,0.2),(-0.4,-0.2),(0.4,0.2)\big\}$, considering a sinusoidal covariate given by $x_t=\sin\big(\pi t/50\big)$ with coefficient $\beta=0.5$. We consider the logit as the link function and a size 100 burn-in is applied in generating the time series. A total of 1,000 replicas of each scenario were generated. Routines to sample from an ULARMA$(p,q)$ process and to perform estimation via PMLE are available in R package \texttt{BTSR}\footnote{The ULARMA model is implemented in the BTSR version 1.0.0 beta, which will be released on CRAN soon.} \citep{btsr}.  For reference, a typical realization of such process is presented in Figure \ref{ref} (top left) along with the respective $\mu_t$.

\begin{figure}[ht!]
\centering
\mbox{
\includegraphics[width=0.98\textwidth]{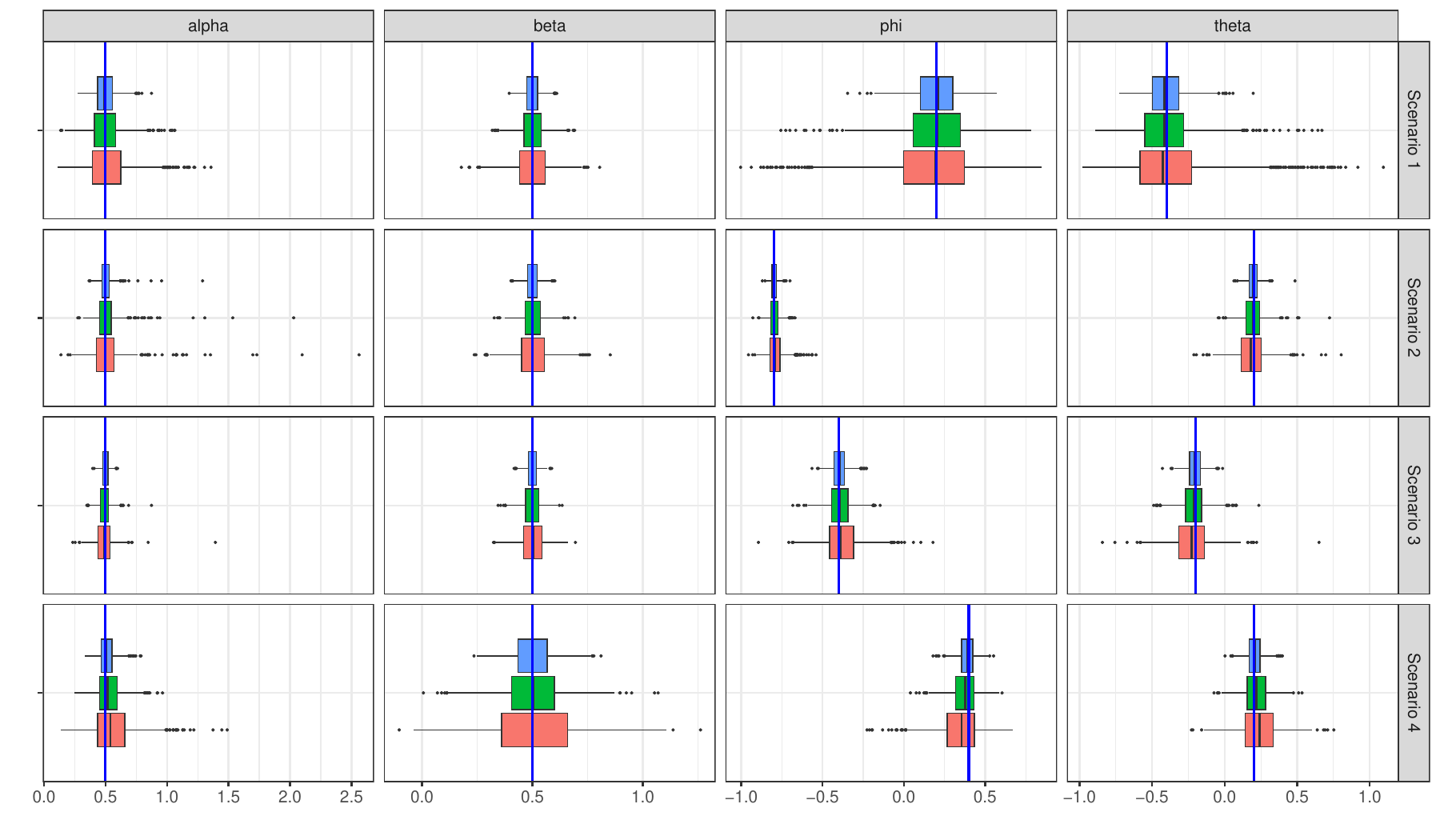}
}
\mbox{
\includegraphics[width=0.98\textwidth]{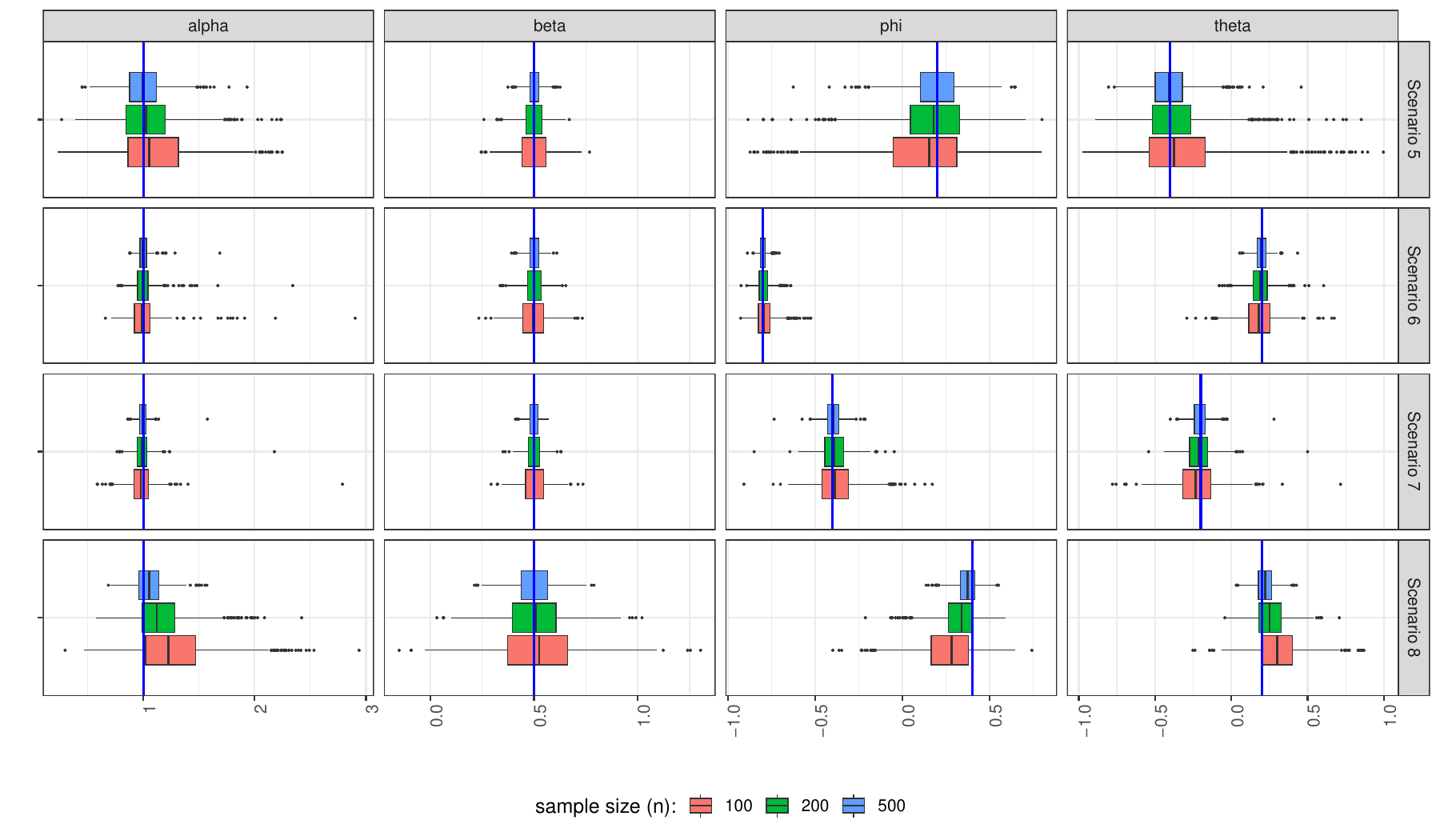}
}
\caption{Boxplots of the simulation results for all parameter for $\alpha=0.5$ (top) and $\alpha=1$ (bottom), with $\beta=0.5$ fixed. Parameter $(\phi,\theta)$ are defined by the scenarios, as follows: scenarios 1 and 5: $(0.2,-0.4)$, scenarios 2 and 6: $(-0.8,0.2)$, scenarios 3 and 7: $(-0.4,-0.2)$, and scenarios 4 and 8: $(0.4,0.2)$.} \label{bp}
\end{figure}

Table \ref{point} summarizes the simulation results. For each set of parameters, we present the mean (left), \emph{median} (center, in italics), and standard deviations (in parentheses) calculated from the 1,000 replicas. Even when $n=100$, parameters $\alpha$ and $\beta$ are well estimated in almost all cases. As for parameters $\phi$ and $\theta$, estimation for $n=100$ is good in some cases and fair in others. For $n=500$, the estimation improves considerably with small biases in most cases. The only exception is $\alpha = 1$, $\phi=0.4$, $\theta=0.2$, which present the worst estimation results among all scenarios. 
As the sample size  $n$ increases, the mean and median estimates tend to converge. This trend is illustrated in Figure \ref{bp}, which presents boxplots of the simulation results. The boxplots corresponding to scenarios where this trend is not observed reveal a significant number of outliers in the estimation, which tend to pull the mean away from the median. Most of these outliers are due to the non-convergence of the L-BFGS-B optimization algorithm. From Table \ref{point} and Figure \ref{bp}, it is evident that the variance of the PMLE vary with the parameter combination and, as expected, diminishes as the sample size increases.
\begin{table}[t!]
\renewcommand{\arraystretch}{1.1}
\setlength{\tabcolsep}{3.6pt}
\caption{Simulation results: point estimates based on 1,000 replicas of each scenario. Presented are the mean (left), the median (center, in italic) and the standard deviation (right, in parentheses).}\label{point}
\centering
\footnotesize
\vspace{.3cm}
\begin{tabular}{c|ccc|ccc|ccc|ccc}
\hline
\multicolumn{1}{c|}{$n$} & \multicolumn{3}{c|}{$\alpha = 0.5$}&\multicolumn{3}{c|}{$\beta = 0.5$} &\multicolumn{3}{c|}{$\phi = 0.2$}&\multicolumn{3}{c}{$\theta = -0.4$}\\
 \hline
100 & 0.521  & \emph{ 0.491 } & $(0.189)$ & 0.498  & \emph{ 0.496 } & $(0.086)$ & 0.158  & \emph{ 0.213 } & $(0.309)$ & -0.371  & \emph{ -0.428 } & $(0.323)$  \\
200 & 0.502  & \emph{ 0.490 } & $(0.142)$ & 0.503  & \emph{ 0.503 } & $(0.059)$ & 0.194  & \emph{ 0.218 } & $(0.225)$ & -0.399  & \emph{ -0.431 } & $(0.225)$  \\
500 & 0.495  & \emph{ 0.490 } & $(0.091)$ & 0.501  & \emph{ 0.502 } & $(0.034)$ & 0.206  & \emph{ 0.215 } & $(0.147)$ & -0.408  & \emph{ -0.419 } & $(0.143)$  \\
\hline
\multicolumn{1}{c|}{$n$} & \multicolumn{3}{c|}{$\alpha = 0.5$}&\multicolumn{3}{c|}{$\beta = -0.5$} &\multicolumn{3}{c|}{$\phi = -0.8$}&\multicolumn{3}{c}{$\theta = 0.2$}\\
 \hline
100 & 0.516  & \emph{ 0.492 } & $(0.190)$ & 0.503  & \emph{ 0.503 } & $(0.074)$ & -0.788  & \emph{ -0.793 } & $(0.052)$ & 0.183  & \emph{ 0.179 } & $(0.109)$  \\
200 & 0.510  & \emph{ 0.499 } & $(0.122)$ & 0.501  & \emph{ 0.499 } & $(0.051)$ & -0.795  & \emph{ -0.797 } & $(0.034)$ & 0.193  & \emph{ 0.190 } & $(0.072)$  \\
500 & 0.505  & \emph{ 0.501 } & $(0.068)$ & 0.501  & \emph{ 0.501 } & $(0.032)$ & -0.798  & \emph{ -0.799 } & $(0.021)$ & 0.197  & \emph{ 0.196 } & $(0.044)$  \\
\hline
\multicolumn{1}{c|}{$n$} & \multicolumn{3}{c|}{$\alpha = 0.5$}&\multicolumn{3}{c|}{$\beta = 0.5$} &\multicolumn{3}{c|}{$\phi = -0.4$}&\multicolumn{3}{c}{$\theta = -0.2$}\\
 \hline
100 & 0.483  & \emph{ 0.483 } & $(0.076)$ & 0.499  & \emph{ 0.497 } & $(0.066)$ & -0.377  & \emph{ -0.386 } & $(0.121)$ & -0.233  & \emph{ -0.232 } & $(0.140)$  \\
200 & 0.495  & \emph{ 0.494 } & $(0.052)$ & 0.499  & \emph{ 0.499 } & $(0.044)$ & -0.391  & \emph{ -0.398 } & $(0.079)$ & -0.213  & \emph{ -0.214 } & $(0.090)$  \\
500 & 0.497  & \emph{ 0.498 } & $(0.031)$ & 0.500  & \emph{ 0.499 } & $(0.027)$ & -0.396  & \emph{ -0.398 } & $(0.047)$ & -0.205  & \emph{ -0.207 } & $(0.055)$  \\
\hline
\multicolumn{1}{c|}{$n$} & \multicolumn{3}{c|}{$\alpha = 0.5$}&\multicolumn{3}{c|}{$\beta = -0.5$} &\multicolumn{3}{c|}{$\phi = 0.4$}&\multicolumn{3}{c}{$\theta = 0.2$}\\
 \hline
100 & 0.578  & \emph{ 0.550 } & $(0.192)$ & 0.512  & \emph{ 0.513 } & $(0.222)$ & 0.329  & \emph{ 0.342 } & $(0.137)$ & 0.248  & \emph{ 0.244 } & $(0.140)$  \\
200 & 0.535  & \emph{ 0.525 } & $(0.117)$ & 0.501  & \emph{ 0.502 } & $(0.151)$ & 0.367  & \emph{ 0.372 } & $(0.085)$ & 0.224  & \emph{ 0.223 } & $(0.092)$  \\
500 & 0.511  & \emph{ 0.508 } & $(0.066)$ & 0.498  & \emph{ 0.501 } & $(0.094)$ & 0.387  & \emph{ 0.389 } & $(0.051)$ & 0.209  & \emph{ 0.208 } & $(0.058)$  \\
\hline
\multicolumn{1}{c|}{$n$} & \multicolumn{3}{c|}{$\alpha = 1$}&\multicolumn{3}{c|}{$\beta = 0.5$} &\multicolumn{3}{c|}{$\phi = 0.2$}&\multicolumn{3}{c}{$\theta = -0.4$}\\
 \hline
100 & 1.086  & \emph{ 1.043 } & $(0.341)$ & 0.502  & \emph{ 0.502 } & $(0.086)$ & 0.127  & \emph{ 0.162 } & $(0.285)$ & -0.341  & \emph{ -0.397 } & $(0.303)$  \\
200 & 1.044  & \emph{ 1.017 } & $(0.263)$ & 0.500  & \emph{ 0.499 } & $(0.060)$ & 0.163  & \emph{ 0.189 } & $(0.218)$ & -0.370  & \emph{ -0.404 } & $(0.220)$  \\
500 & 1.008  & \emph{ 0.999 } & $(0.177)$ & 0.501  & \emph{ 0.500 } & $(0.036)$ & 0.193  & \emph{ 0.199 } & $(0.144)$ & -0.396  & \emph{ -0.406 } & $(0.137)$  \\
\hline
\multicolumn{1}{c|}{$n$} & \multicolumn{3}{c|}{$\alpha = 1$}&\multicolumn{3}{c|}{$\beta = -0.5$} &\multicolumn{3}{c|}{$\phi = -0.8$}&\multicolumn{3}{c}{$\theta = 0.2$}\\
 \hline
100 & 1.002  & \emph{ 0.984 } & $(0.192)$ & 0.503  & \emph{ 0.502 } & $(0.075)$ & -0.784  & \emph{ -0.790 } & $(0.055)$ & 0.177  & \emph{ 0.177 } & $(0.106)$  \\
200 & 1.004  & \emph{ 0.995 } & $(0.133)$ & 0.503  & \emph{ 0.500 } & $(0.053)$ & -0.792  & \emph{ -0.795 } & $(0.035)$ & 0.190  & \emph{ 0.187 } & $(0.072)$  \\
500 & 1.003  & \emph{ 0.998 } & $(0.082)$ & 0.500  & \emph{ 0.501 } & $(0.032)$ & -0.798  & \emph{ -0.799 } & $(0.021)$ & 0.198  & \emph{ 0.199 } & $(0.045)$  \\
\hline
\multicolumn{1}{c|}{$n$} & \multicolumn{3}{c|}{$\alpha = 1$}&\multicolumn{3}{c|}{$\beta = 0.5$} &\multicolumn{3}{c|}{$\phi = -0.4$}&\multicolumn{3}{c}{$\theta = -0.2$}\\
 \hline
100 & 0.980  & \emph{ 0.981 } & $(0.108)$ & 0.501  & \emph{ 0.502 } & $(0.062)$ & -0.383  & \emph{ -0.392 } & $(0.122)$ & -0.223  & \emph{ -0.221 } & $(0.144)$  \\
200 & 0.992  & \emph{ 0.990 } & $(0.072)$ & 0.501  & \emph{ 0.501 } & $(0.043)$ & -0.393  & \emph{ -0.397 } & $(0.081)$ & -0.209  & \emph{ -0.213 } & $(0.093)$  \\
500 & 0.998  & \emph{ 0.998 } & $(0.043)$ & 0.501  & \emph{ 0.502 } & $(0.026)$ & -0.398  & \emph{ -0.400 } & $(0.050)$ & -0.202  & \emph{ -0.203 } & $(0.058)$  \\
\hline
\multicolumn{1}{c|}{$n$} & \multicolumn{3}{c|}{$\alpha = 1$}&\multicolumn{3}{c|}{$\beta = 0.5$} &\multicolumn{3}{c|}{$\phi = 0.4$}&\multicolumn{3}{c}{$\theta = 0.2$}\\
 \hline
100 & 1.297  & \emph{ 1.250 } & $(0.367)$ & 0.513  & \emph{ 0.504 } & $(0.212)$ & 0.259  & \emph{ 0.277 } & $(0.167)$ & 0.307  & \emph{ 0.301 } & $(0.152)$  \\
200 & 1.162  & \emph{ 1.135 } & $(0.242)$ & 0.509  & \emph{ 0.504 } & $(0.149)$ & 0.323  & \emph{ 0.331 } & $(0.111)$ & 0.262  & \emph{ 0.265 } & $(0.107)$  \\
500 & 1.064  & \emph{ 1.062 } & $(0.130)$ & 0.503  & \emph{ 0.503 } & $(0.093)$ & 0.370  & \emph{ 0.372 } & $(0.059)$ & 0.225  & \emph{ 0.226 } & $(0.061)$  \\
\hline
\end{tabular}
\end{table}

\subsection{Goodness-of-fit exercise}
In this section, we evaluate the finite sample performance of the goodness-of-fit tests based on the simple and quantile residuals as discussed in Section \ref{gof}. In addition to the wild bootstrap variance ratio and the Domingu\'ez-Lobato test for the martingale difference hypothesis, we consider five commonly used normality tests for the quantile residuals: the Anderson-Darling (AD), Cram\'er-von Mises (CvM), Kolmogorov-Smirnov (KS), and Shapiro-Francia (SF) tests. To conduct these tests, we apply the R package \texttt{nortest} \citep{nortest} considering all default configurations. For further details on these tests, refer to \cite{thode}.  To perform the wild bootstrap variance ratio test, we generated 500 bootstrap samples using both the Normal and Mammen's two-point distributions (WB). For details, please refer to the documentation of the \texttt{vrtest} package. The Dom\'inguez-Lobato test is also implemented in the \texttt{vrtest} package. However, the existing implementation is written purely in R, making it too slow for simulations purposes. To address this limitation, we re-implemented the \texttt{DL.test} function from \texttt{vrtest} in FORTRAN, which is then called from R to perform the test. The code runs over 90 times faster than the original and is available upon request. In the simulation, we consider two approaches: the Cram\'er-von Mises (DL-Cp) and the Kolmogorov-Smirnov (DL-KS) tests.

We generate samples of size $n\in\{100,200,500\}$ from an ULARMA$(1,1)$ model with parameters $\alpha\in\{0.5,1\}$, $(\phi,\theta) \in\big\{(0.2,-0.8),(-0.8,0.2),(-0.4,-0.2),(0.4,0.2)\big\}$. A size 100 burn-in is applied in generating the time series. The logit link function is used in our experiments. A total of 1,000 replicas of each scenario were generated and tests were evaluated at 5\% confidence level.  
Table \ref{tests} provides a summary of the simulation results. All normality tests demonstrated satisfactory performance, with rejection rates close to the nominal value of 0.05 for most parameters. The only notable exception occurred in the scenario where $\phi=-0.8$ and $\theta=0.2$, where the SF test exhibited slightly higher-than-expected rejection rates. The martingale distance tests presented a conservative performance, with the DL variants presenting rejection rates slightly closer to 0.05 than the wild bootstrap. The results using the normal probability option in the function vrtest were very similar to the ones using the Mammen's two point distribution and were omitted.
\begin{table}[ht!]
\renewcommand{\arraystretch}{1.2}
\setlength{\tabcolsep}{4.8pt}
\caption{Simulation results - martingale difference and normality tests. Presented are the proportion of tests that rejected the null hypothesis for each specific test. The tests were performed at level $\delta=0.05$.}\label{tests}
\centering
\tiny
\vspace{.3cm}
\begin{tabular}{c|c||c|c|c||c|c|c|c||c|c|c||c|c|c|c}                                                
\hline                                                                                             
\hline
  \multirow{2}{*}{$n$} &\multirow[c]{2}{*}{$(\phi,\theta)$}&\multicolumn{7}{c||}{$\alpha = 0.5$}&\multicolumn{7}{c}{$\alpha = 1$}\\
  \cline{3-16} 
  & & WB & DL-Cp & DL-Kp & AD & CvM & KS & SF & WB & DL-Cp & DL-Kp & AD & CvM & KS & SF \\         
  \hline
100 & \multirow[c]{3}{1.5cm}{\begin{tabular}{l}
        $\phi = 0.2$\\ $\theta = -0.4$ \\
        \end{tabular}}
       & 0.00  & 0.01 &  0.02 &   0.04 &  0.05 &  0.04 &  0.05 & 0.01 & 0.02 &  0.03 &    0.06 &  0.06 &  0.05 & 0.06 \\ 
200 &  & 0.00  & 0.00 &  0.01 &   0.04 &  0.05 &  0.05 &  0.05 & 0.01 & 0.02 &  0.02 &    0.05 &  0.05 &  0.05 & 0.04 \\ 
500 &  & 0.00  & 0.01 &  0.02 &   0.06 &  0.06 &  0.07 &  0.05 & 0.01 & 0.01 &  0.01 &    0.05 &  0.05 &  0.05 & 0.05 \\ 
\hline
100 & \multirow[c]{3}{1.5cm}{\begin{tabular}{l}
        $\phi = -0.8$\\ $\theta = 0.2$ \\
        \end{tabular}}
       & 0.02  & 0.03 &  0.04 &   0.09 &  0.08 &  0.07 &  0.13 & 0.03 & 0.02 &  0.04 &    0.09 &  0.08 &  0.06 & 0.10 \\
200 &  & 0.02  & 0.03 &  0.03 &   0.08 &  0.07 &  0.06 &  0.11 & 0.02 & 0.02 &  0.03 &    0.07 &  0.06 &  0.06 & 0.11 \\
500 &  & 0.02  & 0.02 &  0.03 &   0.05 &  0.05 &  0.05 &  0.11 & 0.02 & 0.02 &  0.03 &    0.07 &  0.07 &  0.06 & 0.09 \\
\hline
100 & \multirow[c]{3}{1.5cm}{\begin{tabular}{l}
        $\phi = -0.4$\\ $\theta = -0.2$ \\
        \end{tabular}}
       & 0.01  & 0.02 &  0.03 &   0.07 &  0.07 &  0.06 &  0.08 & 0.01 & 0.01 &  0.02 &    0.06 &  0.06 &  0.06 & 0.06 \\
200 &  & 0.01  & 0.01 &  0.02 &   0.06 &  0.05 &  0.05 &  0.07 & 0.01 & 0.01 &  0.02 &    0.06 &  0.06 &  0.04 & 0.06 \\
500 &  & 0.01  & 0.01 &  0.02 &   0.05 &  0.05 &  0.04 &  0.06 & 0.01 & 0.01 &  0.02 &    0.05 &  0.06 &  0.04 & 0.06 \\
\hline
100 & \multirow[c]{3}{1.5cm}{\begin{tabular}{l}
        $\phi = 0.4$\\ $\theta = 0.2$ \\
        \end{tabular}}
       & 0.01  & 0.02 &  0.02 &   0.05 &  0.05 &  0.06 &  0.07 & 0.05 & 0.03 &  0.03 &    0.05 &  0.06 &  0.05 & 0.05 \\  
200 &  & 0.01  & 0.01 &  0.01 &   0.06 &  0.06 &  0.05 &  0.07 & 0.04 & 0.02 &  0.04 &    0.06 &  0.06 &  0.06 & 0.09 \\  
500 &  & 0.00  & 0.01 &  0.01 &   0.05 &  0.05 &  0.05 &  0.07 & 0.02 & 0.01 &  0.02 &    0.06 &  0.06 &  0.05 & 0.13 \\  
\hline
\end{tabular}
\end{table}

\subsection{Joint behavior}
\begin{figure}[t!]
\centering
\mbox{
\includegraphics[width=.45\textwidth]{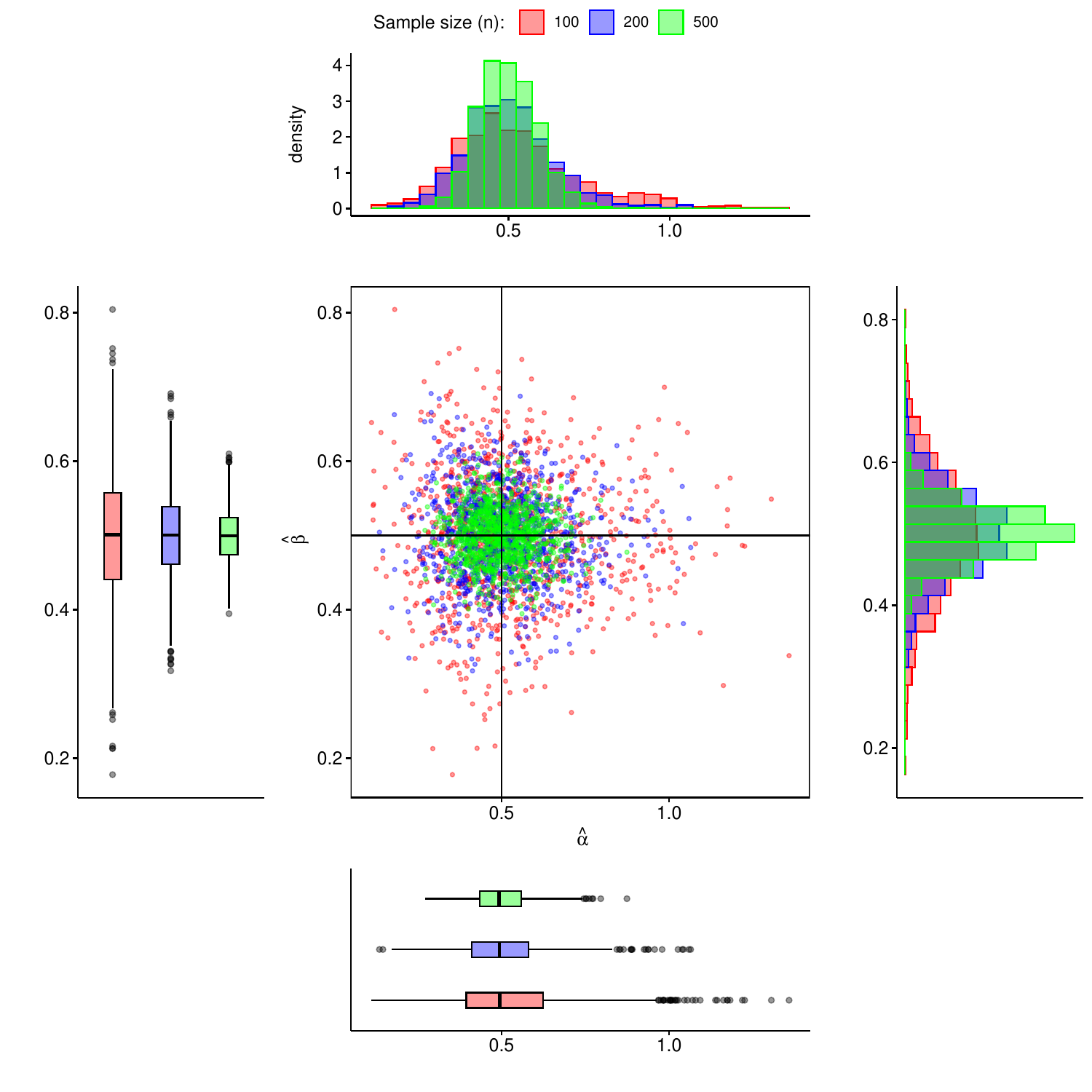}
\includegraphics[width=.45\textwidth]{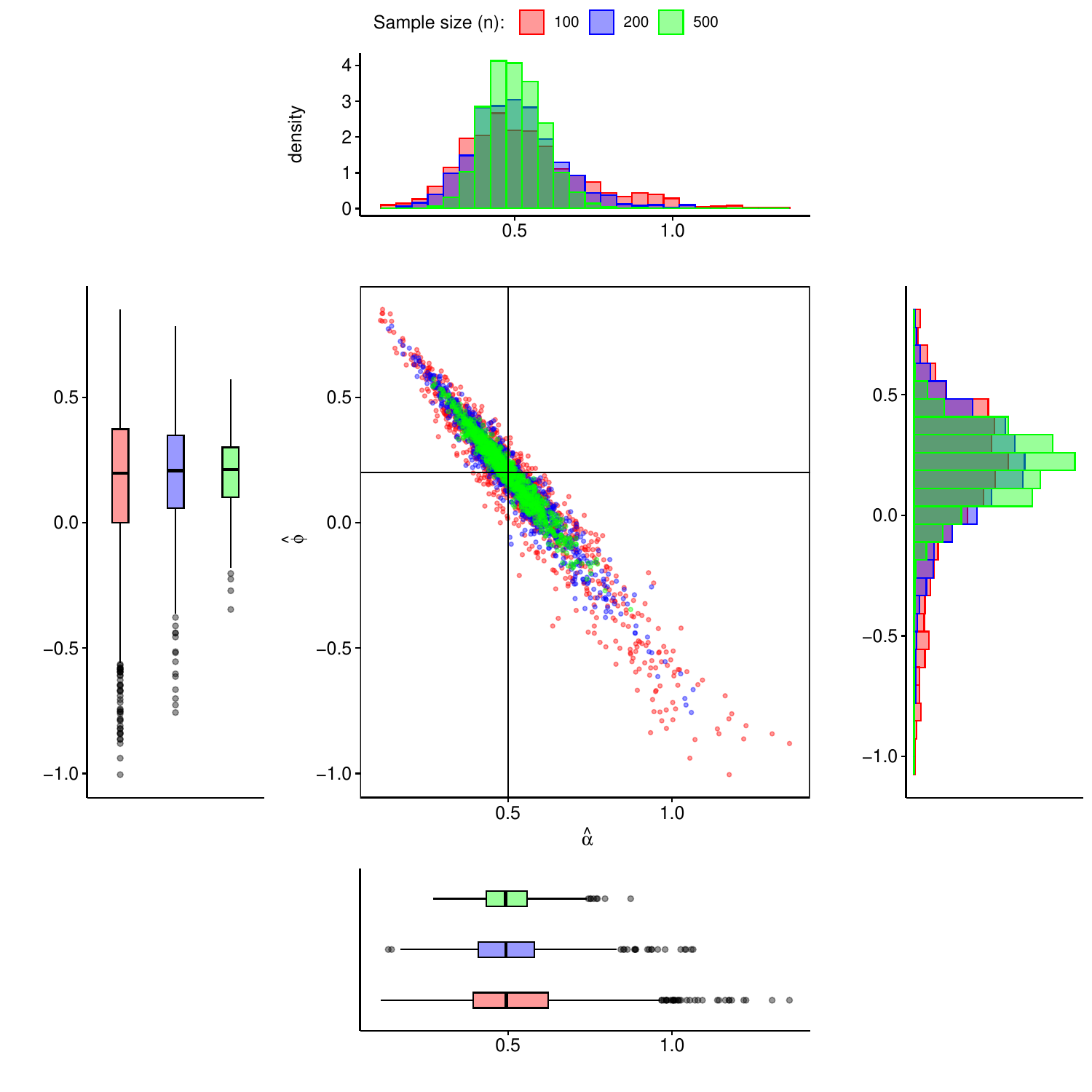}
}
\mbox{
\includegraphics[width=.45\textwidth]{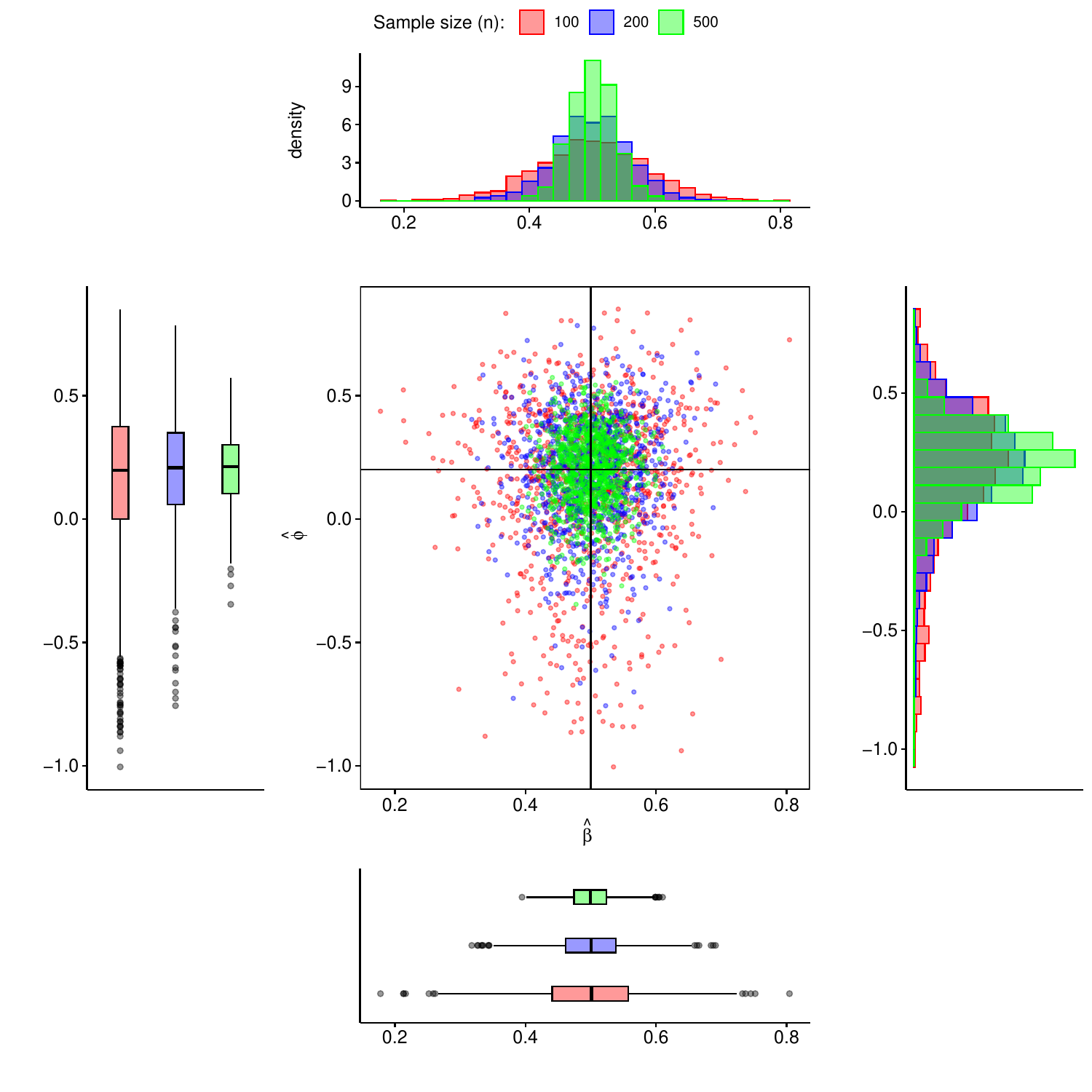}
\includegraphics[width=.45\textwidth]{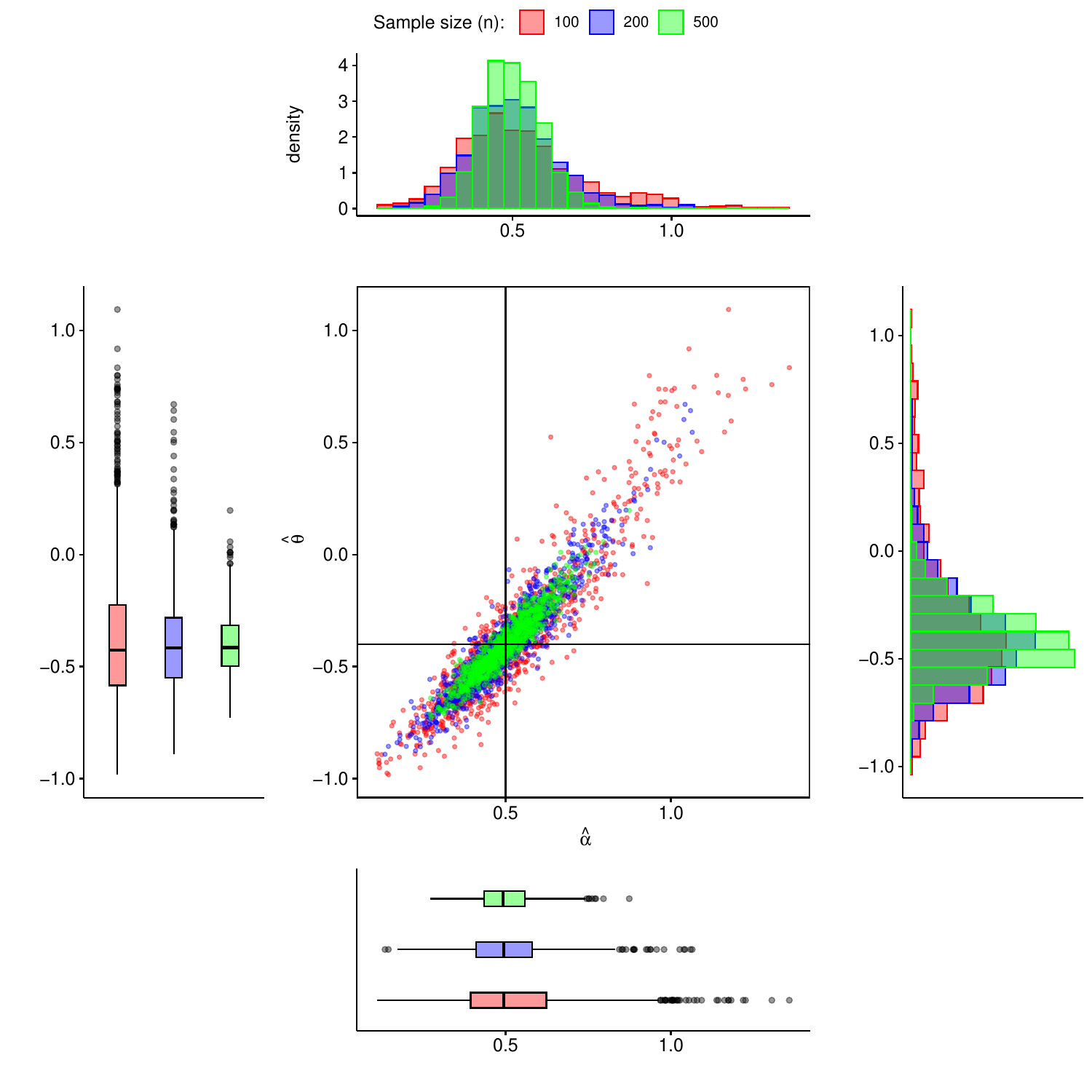}
}
\mbox{
\includegraphics[width=.45\textwidth]{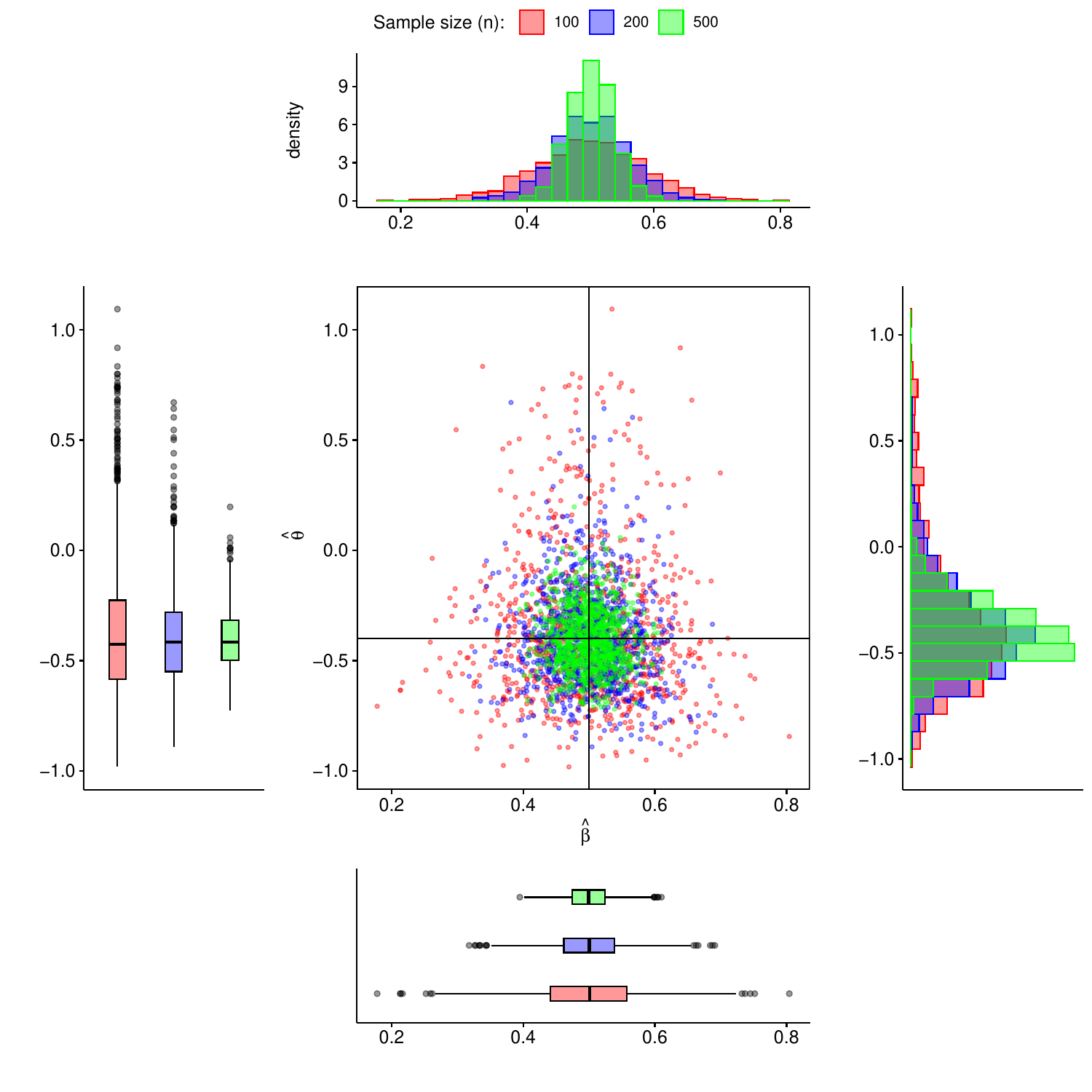}
\includegraphics[width=.45\textwidth]{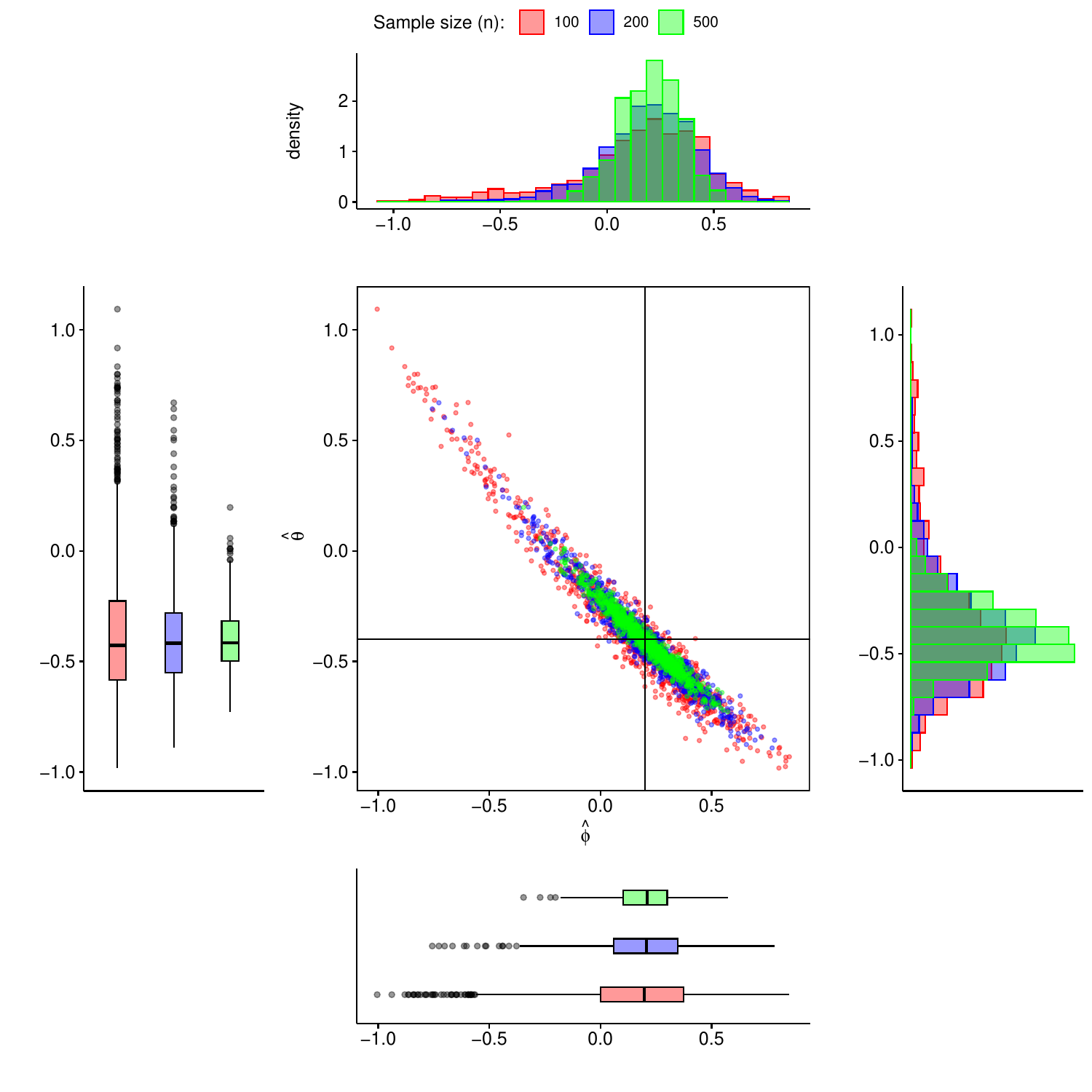}
}
\caption{Pairwise joint and marginal behavior of the estimated values for $\alpha=0.5$, $\beta=0.5$, $\phi=0.2$, $\theta=-0.4$. Solid lines in the scatter plot represent the true values.}
\label{pair}
\end{figure}

The asymptotic normality of the PMLE, as developed in Section \ref{lsi}, is also investigated. To achieve this, we examine pairwise scatter plots and marginal behavior (histograms and boxplots) shown in Figure \ref{pair} for the case $\alpha=0.5$, $\beta=0.5$, $\phi=0.2$, $\theta=-0.4$. Similar behavior is observed in other cases. From the scatter plots the convergence to a bivariate Gaussian pattern as the sample size $n$ increases is clear. Additionally, the marginal behavior (as depicted in the histograms and boxplots) becomes more symmetrical with increasing $n$, and the histograms increasingly resemble the shape of a normal distribution. It is also noteworthy that there is a dependence between the estimates of the parameters $\alpha$, $\phi$, and $\theta$, as evidenced by the pairwise plots located in the right column of Figure \ref{pair}. In stark contrast, the estimates of $\beta$ appear uncorrelated with the other parameters, as shown in the plots in the left column of Figure \ref{pair}.
\FloatBarrier
\section{Application to Real Data}
In this section, we evaluate and compare the performance of the proposed ULARMA model in forecasting the proportion of net electricity generated by conventional hydroelectric power in the United States. Such analysis is important for multiple reasons. First, hydroelectric power generation fluctuates due to climatic and hydrological variability, making accurate forecasting essential for mitigating uncertainties and optimizing energy distribution. Second, as the energy sector increasingly shifts toward renewable sources, understanding hydroelectric generation patterns and its participation in net power generation allows for a more effective energy resources management. Finally, by benchmarking the proposed method against existing models, we can assess its predictive accuracy, reliability, and practical applicability in real-world scenarios. %This comparison not only validates the effectiveness of the proposed ULARMA models but also highlights its potential advantages.

\subsubsection*{Data description}
The dataset used in this section represent the monthly proportion of net electricity generated by conventional hydroelectric for all sectors in United States. The raw data is freely available from the U.S. Energy Information Administration's (EIA) website\footnote{{\color{blue} \href{https://www.eia.gov/electricity/data/browser/\#/topic/}{www.eia.gov/electricity/data/browser/\#/topic/}} retrieved 04/04/2025.} and contains monthly net energy generated for all sector of United States considering coal, petroleum liquids, petroleum coke, natural gas, other gases, nuclear, conventional hydroelectric, other renewables (total), hydro-electric pumped storage and other. Our interested lies on the ratio between the net electricity generated by hydroelectric and the total net energy generated by all fuels considered in the data from January 2001 to January 2025. For forecasting purposes, one year data from February 2024 to January 2025 was reserved. Hence, only the data from January 2001 to January 2024 was used for fitting purposes, yielding a time series of size $n=277$. 

\begin{figure}[!ht]
\centering
\subfigure[Observed data]{\includegraphics[width=0.49\textwidth]{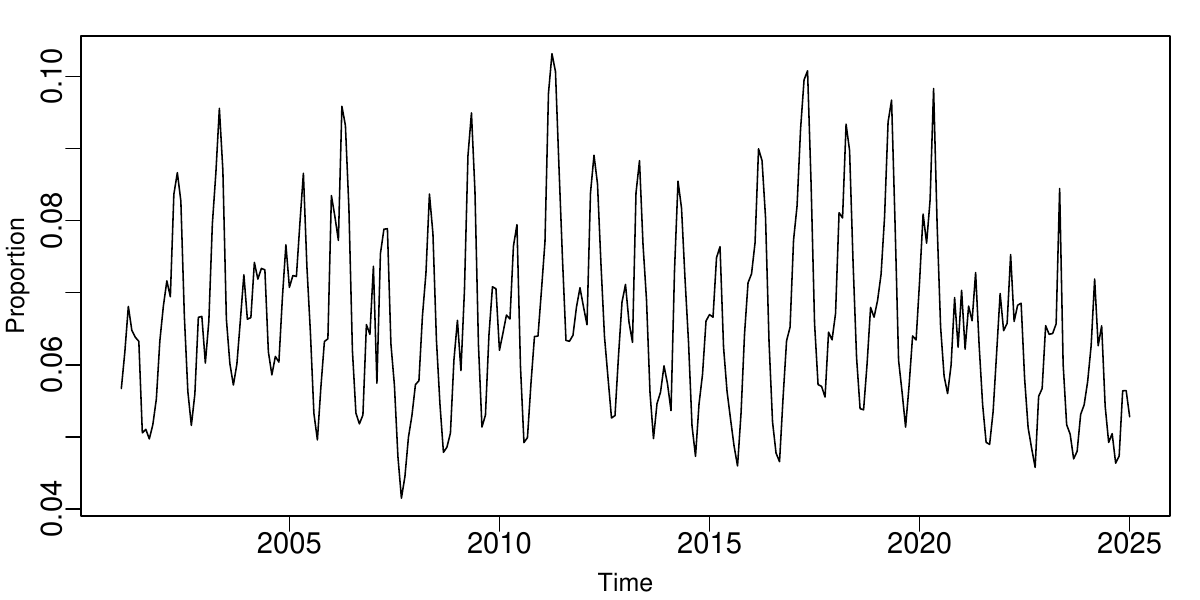}}
\subfigure[Sample ACF]{\includegraphics[width=0.49\textwidth]{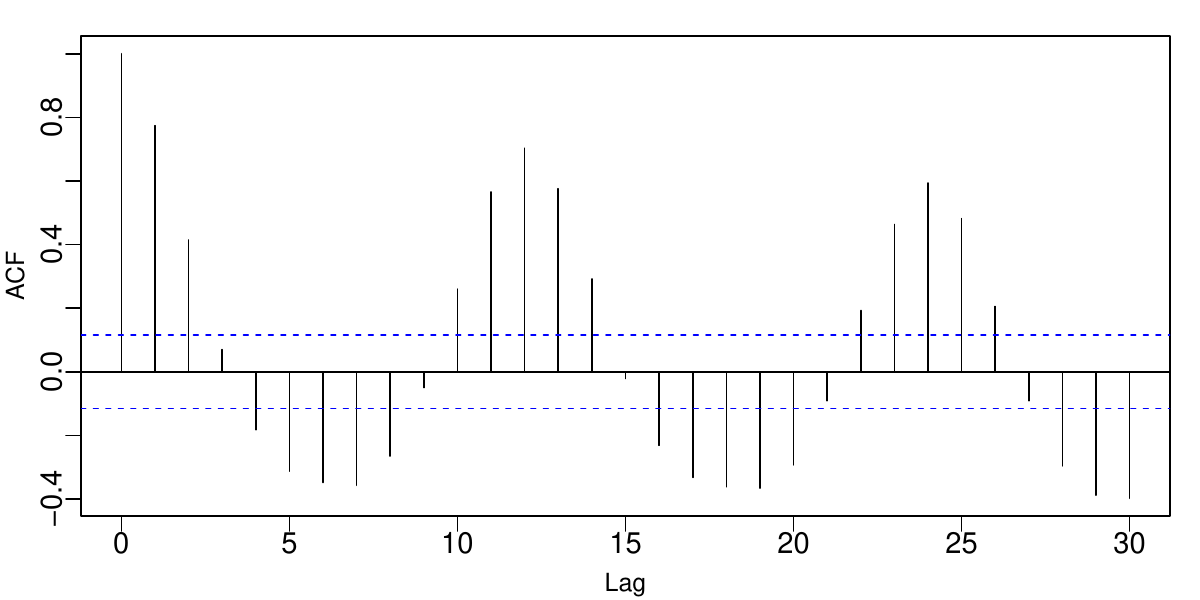}}
\subfigure[Sample PACF]{\includegraphics[width=0.49\textwidth]{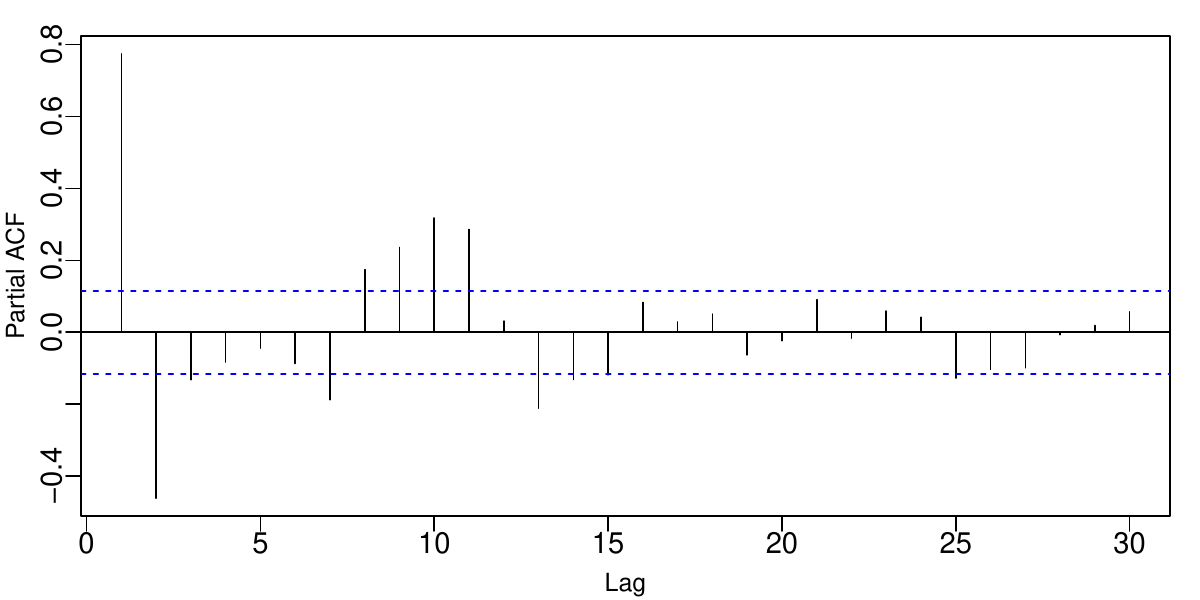}}
\subfigure[Seasonality]{\includegraphics[width=0.49\textwidth]{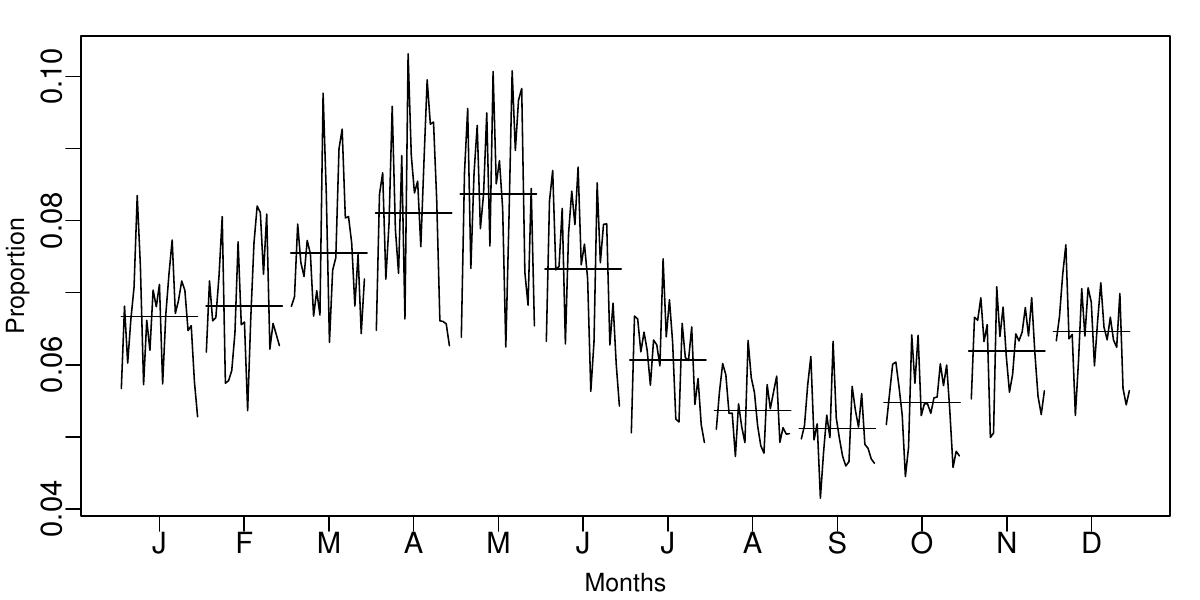}}
\caption{Proportion of net energy generated by hydroelectric data: (a) the observed time series and the corresponding (b) autocorrelation function (ACF), (c) partial autocorrelation (PACF) function and (d) seasonal plot.}\label{fig:f}
\end{figure}

The time series plot of the complete data is presented in Figure \ref{fig:f}. A seasonal cycle in the proportion of energy produced by hydroelectric power generation is clearly discernible in the plots. This seasonality is primarily driven by hydrological and climatic factors, as well as energy demand variations throughout the year. There are two main approaches to model time series presenting seasonality. The first one is by introducing deterministic covariates combining sine and cosine to model the year seasonal pattern, whereas the second one is to model the data directly without the aid of covariates. In this work we are going to follow the second approach to highlight the ULARMA's versatility in modeling seasonal patterns without the aid of seasonal components and its superior forecasting capabilities in this scenario. To further highlight that, we shall compare the ULARMA to benchmark models -- the KARMA and $\beta$ARMA models -- which are based on bi-parametric conditional distribution and are known to be very flexible. 
\FloatBarrier

\subsection{Parameter Estimation, Inference and Model Identification}\label{id}

Parameter estimation was performed using the L-BFGS-B and Nelder-Mead algorithms combined, as implemented in the \texttt{BTSR} package. We tested the logit, log-log, and clog-log link functions to determine the most suitable one. Model selection followed a bidirectional stepwise procedure. The process began with a full ULARMA$(12,12)$ model. During the backward elimination phase, parameters with $p$-values exceeding 15\% were sequentially fixed at zero, ensuring that only statistically significant terms were retained. In the forward selection phase, previously excluded parameters were reassessed and reintroduced if their $p$-values fell below the 10\% significance threshold. The choice of two different thresholds for inclusion and exclusion of parameter was chosen to avoid the procedure to being stuck removing and adding the same parameter in a loop. This iterative process, guided by the Wald test described in Section~\ref{s:hp}, continued until no further improvements were identified. 

Additionally, we choose two benchmark models - the KARMA and $\beta$ARMA models - to compare with the ULARMA. The $\beta$ARMA model was one of the first GARMA-like models introduced in the literature, whereas the KARMA is deemed a more robust alternative to the $\beta$ARMA, being popular in the hydrological literature. Both are very flexible models based on bi-parametric conditional distributions and are part of the \texttt{BTSR} package. For these benchmarks, the same bidirectional stepwise procedure was applied for model selection. 

To provide a comparison of the predicted values among the competing models, we computed some in-sample goodness-of-fit measures, such as root mean squared error~(RMSE), mean absolute percentage error~(MAPE), and mean directional accuracy~(MDA). The MDA evaluates the model's ability to correctly predict the direction of change in the time series, which is particularly relevant in applications where capturing trend reversals and directional movements is important. Higher values of the MDA reflect improved alignment of the forecasts with upward or downward movements of the time series.

%For each model fitted, we chose the link for which the RMSE was smaller for the particular model.
\begin{table}[ht]
\centering
\caption{Summary results for the best model related to each link function: order ($p$ and $q$), number of significant parameters ($k$), $p$-values from DL and SF tests, in-sample goodness-of-fit measures (RMSE, MAPE, MDA) and AIC.}\vskip 0.5\baselineskip
\label{tab:results}
\setlength{\tabcolsep}{4pt}
\renewcommand{\arraystretch}{1.4}
\begin{tabular}{c|c|ccc|cc|ccc|c}
Model & Link & $p$ & $q$ & $k$  &  DL & SF & RMSE & MAPE & MDA & AIC  \\
\hline   
\multirow{3}{*}{ULARMA}
   & logit   & 11 & 0 & 4     & 0.3357 & 0.5243         & 0.0101 & 0.0984 & 0.6558   & -982.6\\
   & loglog  & 4  & 0 & 3     & 0.3502 & 0.0813         & 0.0097 & 0.0974 & 0.6522   & -984.5\\
   & cloglog & 11 & 0 & 5     & 0.0686 & 0.7044         & 0.0082 & 0.0909 & 0.6920   & -981.8\\
\hline                                                    
\multirow{3}{*}{ KARMA}                                   
   & logit   & 11 & 0 & 7     & 0.4874 &$<10^{-5}$& 0.0094 & 0.0819 & 0.7790   & -1977.9\\
   & loglog  & 11 & 0 & 7     & 0.3610 & $<10^{-5}$& 0.0084 & 0.0835 & 0.7717   & -1970.7\\
   & cloglog & 11 & 0 & 7     & 0.5343 & $<10^{-5}$& 0.0126 & 0.0834 & 0.7681   & -1982.1\\
\hline                                                    
\multirow{3}{*}{$\beta$ARMA}                              
   & logit   & 9  & 12 & 11   & 0.3791 & 0.0250         & 0.0068 & 0.0792 & 0.7572   & -1854.1\\
   & loglog  & 11 & 1  & 10   & 0.3105 & 0.0805         & 0.0068 & 0.0825 & 0.7391   & -1894.6\\
   & cloglog & 8  & 12 & 11   & 0.3971 & 0.0179         & 0.0069 & 0.0803 & 0.7536   & -1836.5\\

\end{tabular}
\end{table}

Simulation resultsare presented in Table \ref{tab:results}, which reveals some interesting results. Regarding model structure, 6 out of the 9 models considered $p=11$ as the ``leading'' autoregressive order, whereas only the $\beta$ARMA models included moving average terms.  The ULARMA were the most parsimonious models including between 3 and 5 significant coefficients in the model. The KARMA models followed in second, with exactly 7 coefficients for each model, whereas the $\beta$ARMA presented uniformly more complex models, with 10 to 11 significant coefficients in each model. 

Considering simple residuals, none of the models rejected the martingale difference null hypothesis considering the Dom\'inguez-Lobato test. It is interesting to notice that even the simple residuals for the KARMA, which theoretically have no reason to follow a martingale difference, as $\mu_t$ represents the median in this case, did not reject the null hypothesis. We also performed the Shapiro-Francia (SF) normality test on the  quantile residuals for each fitted model. The quantile residuals of all models presented a very pronounced outlier at $t=1$ due to initialization, so that the SF test was applied removing the first observation from the quantile residuals for all models. For the ULARMA models, the quantile residuals of all model configurations did not reject the null hypothesis of normality. The same happens for the $\beta$ARMA models based on log-log link function, whereas the other configurations rejected the null hypothesis of the SF test. The quantile residuals for the KARMA model presented significative departure from normality, which yielded very low $p-$values for the SF test.

Regarding in-sample goodness-of-fit, the $\beta$ARMA presented uniformly better results in terms of RMSE and MAPE, whereas the ULARMA presented the worst. This may be a consequence of model parsimony.  As for the MDA, KARMA models presented the uniformly higher values, hence being the best performing model in MDA sense. Again, ULARMA presented the worst result. AIC within each model are close to each other being hard to advocate for the use of one or other link function solely based on the AIC. Finally, considering the effects of the link in terms of in-sample goodness-of-fit (and AIC), the results were close to each other within the same model, suggesting no influence of the link in this regard. 

\begin{table}[!ht]
\centering
\caption{Fitted ULARMA, KARMA and $\beta$ARMA models: estimated coefficients and respective standard errors (in parenthesis).}\label{t:cacondefit}\vskip 0.5\baselineskip
\renewcommand{\arraystretch}{1.1}
\small
\begin{tabular}{c|ccc|ccc|ccc}
\hline
Model &  \multicolumn{3}{c}{ULARMA} & \multicolumn{3}{|c}{KARMA}& \multicolumn{3}{|c}{$\beta$ARMA}\\
\cline{2-10}
Link &  logit & log-log & clog-log &  logit & log-log & clog-log&  logit & log-log & clog-log\\
\hline
\multirow{2}{*}{ $\alpha$ } 
 &  -1.708 &  0.758 &  -2.390 &  -1.573 &  0.651 &  -1.297 &  -2.384 &  1.025 &  -2.449 \\ 
 & (0.633) & (0.243) & (0.646) & (0.072) & (0.029) & (0.071) & (0.133) & (0.047) & (0.133) \\ 
 \hline
\multirow{2}{*}{ $\phi_{1}$ } 
 &  0.433 &  0.618 &  0.482 &  0.980 &  0.741 &  0.918 &  0.726 &  0.599 &  0.748 \\ 
 & (0.219) & (0.228) & (0.229) & (0.034) & (0.034) & (0.033) & (0.077) & (0.051) & (0.073) \\ 
 \hline
\multirow{2}{*}{ $\phi_{2}$ } 
 &   -  &   -  &   -  &  -0.455 &  -0.207 &  -0.447 &  -0.182 &   -  &  -0.183 \\ 
&  - &  - &  -  & (0.034) & (0.032) & (0.028) & (0.070)&  -  & (0.068) \\ 
\hline
\multirow{2}{*}{ $\phi_{3}$ } 
 &   -  &   -  &   -  &   -  &   -  &   -  &   -  &  -0.176 &   -  \\ 
&  - &  - &  - &  - &  - &  - &  -  & (0.053)&  -  \\ 
\hline
\multirow{2}{*}{ $\phi_{4}$ } 
 &   -  &  -0.356 &  -0.305 &  -0.086 &  -0.138 &   -  &  -0.179 &  -0.189 &  -0.223 \\ 
&  -  & (0.146) & (0.177) & (0.023) & (0.022)&  -  & (0.050) & (0.049) & (0.050) \\ 
\hline
\multirow{2}{*}{ $\phi_{6}$ } 
 &   -  &   -  &   -  &   -  &   -  &   -  &  0.141 &   -  &  0.100 \\ 
&  - &  - &  - &  - &  - &  -  & (0.058)&  -  & (0.058) \\ 
\hline
\multirow{2}{*}{ $\phi_{7}$ } 
 &  -0.299 &   -  &   -  &  -0.196 &  -0.245 &  -0.172 &  -0.327 &  -0.154 &  -0.286 \\ 
 & (0.158)&  - &  -  & (0.019) & (0.018) & (0.021) & (0.060) & (0.050) & (0.066) \\
 \hline
\multirow{2}{*}{ $\phi_{8}$ } 
 &   -  &   -  &  -0.280 &   -  &   -  &   -  &   -  &  -0.109 &  -0.080 \\ 
&  - &  -  & (0.180)&  - &  - &  - &  -  & (0.059) & (0.053) \\ 
\hline
\multirow{2}{*}{ $\phi_{9}$ } 
 &   -  &   -  &   -  &   -  &   -  &  -0.104 &  -0.087 &  -0.078 &   -  \\ 
&  - &  - &  - &  - &  -  & (0.025) & (0.041) & (0.054)&  -  \\ 
\hline
\multirow{2}{*}{ $\phi_{11}$ } 
 &  0.244 &   -  &  0.234 &  0.161 &  0.194 &  0.319 &   -  &  0.078 &   -  \\ 
 & (0.133)&  -  & (0.148) & (0.014) & (0.013) & (0.019)&  -  & (0.035)&  -  \\ 
 \hline
\multirow{2}{*}{ $\theta_{1}$ } 
 &   -  &   -  &   -  &   -  &   -  &   -  &  0.282 &  0.234 &  0.330 \\ 
&  - &  - &  - &  - &  - &  -  & (0.116) & (0.105) & (0.108) \\ 
\hline
\multirow{2}{*}{ $\theta_{11}$ } 
 &   -  &   -  &   -  &   -  &   -  &   -  &  0.325 &   -  &  0.338 \\ 
&  - &  - &  - &  - &  - &  -  & (0.093)&  -  & (0.091) \\
\hline
\multirow{2}{*}{ $\theta_{12}$ } 
 &   -  &   -  &   -  &   -  &   -  &   -  &  0.401 &   -  &  0.429 \\ 
&  - &  - &  - &  - &  - &  -  & (0.097)&  -  & (0.095) \\ 
\hline
\multirow{2}{*}{ $\nu$ } 
 &   -  &   -  &   -  &  11.054 &  11.012 &  11.028 &  467.98 &  591.39 &  431.18 \\ 
&  - &  - &  -  & (0.517) & (0.515) & (0.516) & (39.807) & (50.297) & (36.678) \\ 
\hline
\end{tabular}
\end{table}
From Table \ref{t:cacondefit}, it is interesting to notice that AR and MA coefficients present in the final selected model present the same signal across all models. For $\alpha$, the estimated values obtained across all models are always negative for  the logit and clog-log, and positive for log-log. The only parameter (besides $\alpha$) present in all models is $\phi_1$. The second most commonly significative parameter was $\phi_4$ and $\phi_7$ present in 7 out of the 9 models. Parameter $\theta$ was significant for all fitted $\beta$ARMA models, but for no other.

\FloatBarrier
\subsection{Forecasting Exercise}

\begin{figure}[!ht]
\centering
\subfigure[In-sample forecast]{\includegraphics[width=0.47\textwidth]{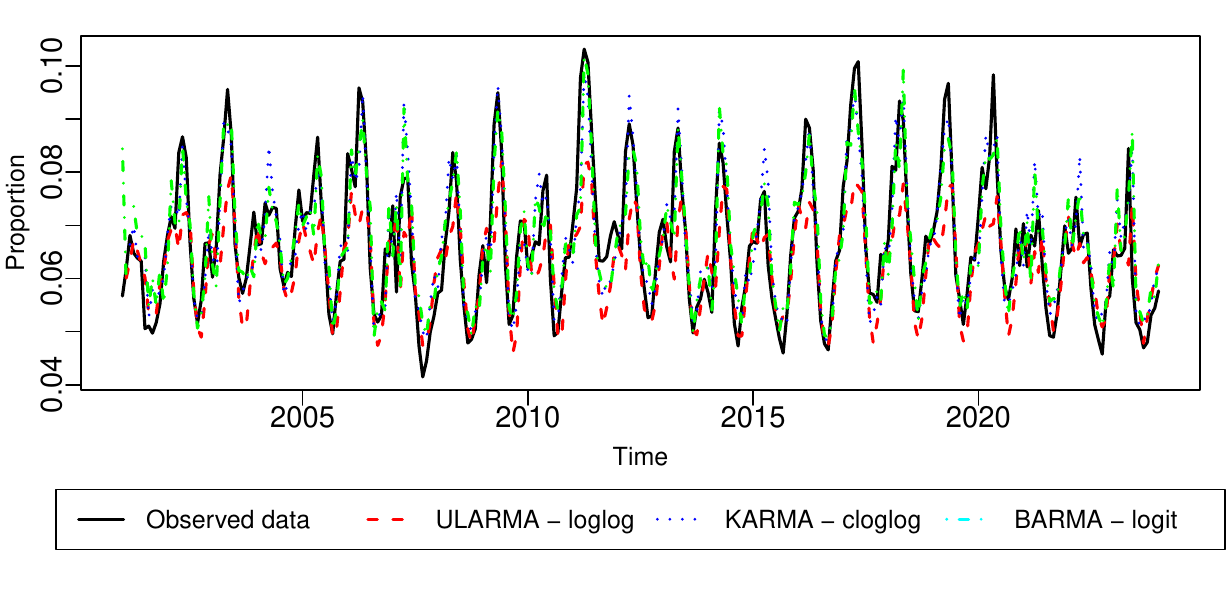}
\label{f:fit}}
\subfigure[Out-of-sample forecast]{\includegraphics[width=0.47\textwidth]{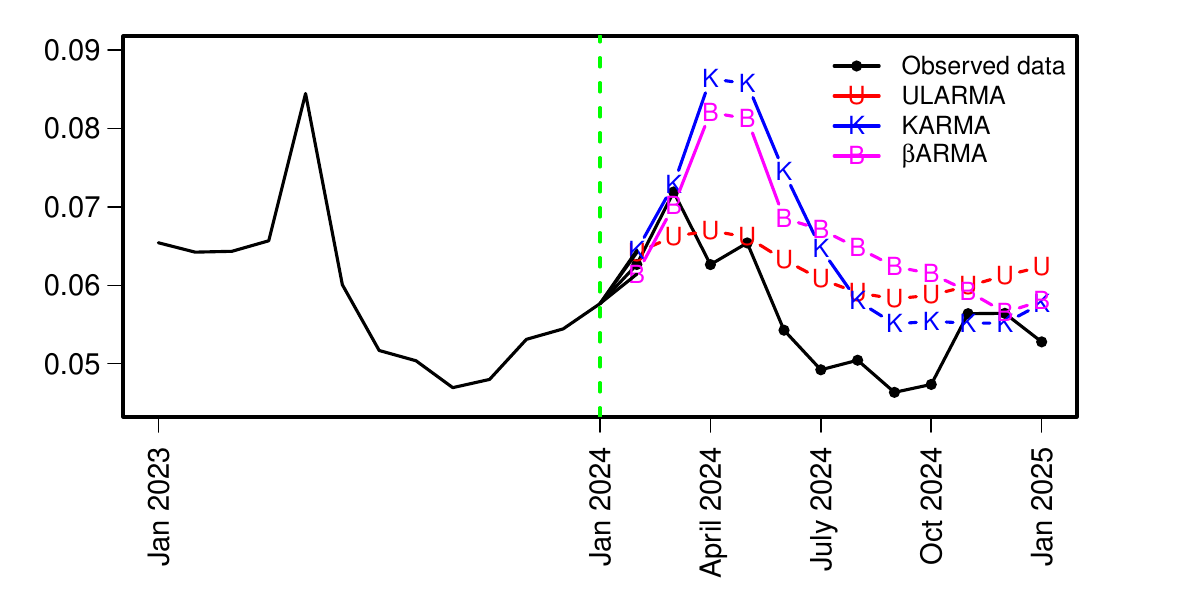}
 \label{f:prev}}
 \subfigure[Bootstrap PI]{\includegraphics[width=0.47\textwidth]{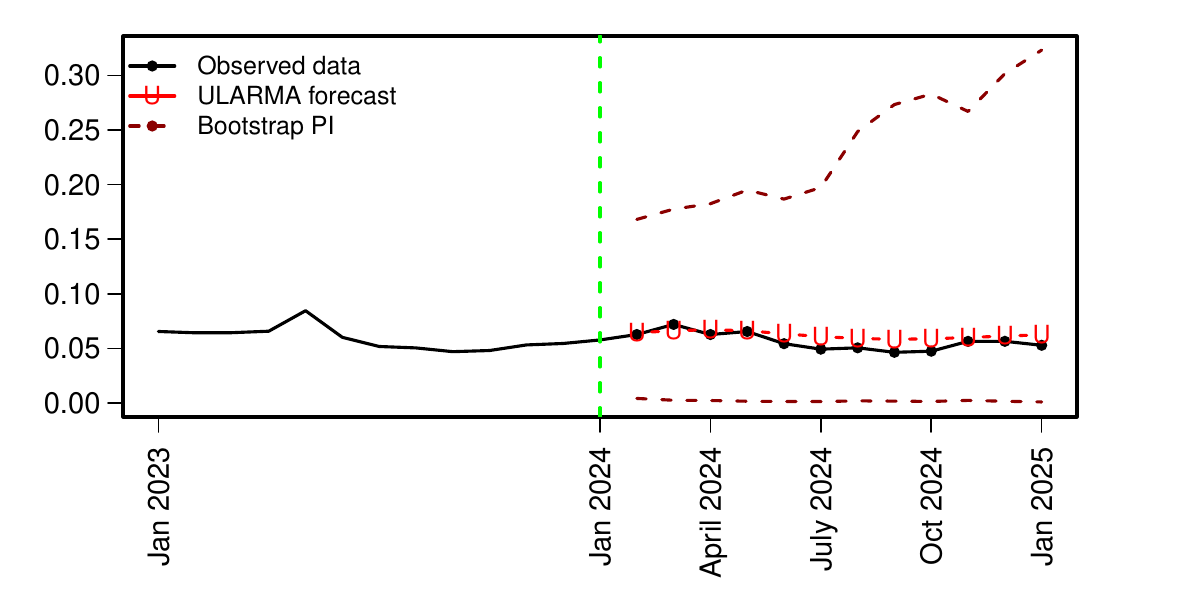}
 \label{f:prevB}}
\caption{Observed time series, (a) fitted values (in-sample-forecast) and (b) 12-step ahead (out-of-sample) forecasts obtained with the fitted ULARMA, KARMA and $\beta$ARMA models. (c) Level 10\% Bootstrap prediction intervals.}\label{f:adjusted}
\end{figure}

Forecasting is probably the most important goal of time series analysis. In this section we consider 12-steps ahead forecast of the monthly proportion of net electricity generated by conventional hydroelectric for all sectors in United States. For simplicity, among the 9 fitted models presented in Section \ref{id}, we shall conduct the analysis considering only the best predicting model of each type. These are the ULARMA considering the log-log, KARMA with the clog-log and $\beta$ARMA with the logit link. The fitted models are presented in Table 4. The coefficients of the three models presented $p-$values smaller than 0.035. Figure \ref{f:adjusted} presents the observed time series, fitted values (in-sample forecast)  and 12-steps ahead (out-of-sample) forecasts. From Figure \ref{f:fit}, we observe that the in-sample forecasts were generally close to the observed values. The fitted values for the KARMA and $\beta$ARMA are difficult to distinguish in the plots, whereas the ULARMA's are easier to spot. 

As for the out-of-sample forecasts, from Figure \ref{f:prev} we observe that for the first two steps the models produce somewhat comparable predictions but after the third it is clear that the KARMA and $\beta$ARMA seem to predict a peak in the data that doesn't exist, which may be an indication of overfitting. ULARMA's forecasted values stayed at lower values, much closer to the observed data. The KARMA and $\beta$ARMA approach the observed data in the end of the 12-step horizon. 

Table \ref{t:f} complements the results, presenting the RMSE, MAPE forecast accuracy measures. Considering the RMSE and MAPE, the ULARMA present the best performance for all except the first two steps, where the $\beta$ARMA and KARMA divide the lead. Considering 12-step ahead forecasts, the ULARMA presented a MAPE of 13.18\%, followed by the KARMA with 17.18\% and the $\beta$ARMA with 19.2\%. Figure \ref{f:prevB} present the level 10\% bootstrap prediction interval obtained using the methodology presented in Section \ref{spi}, calculated from 500 bootstrap samples. We observe that the upper confidence bounds are considerably large, probably due to the presence of several peaks on the data, a feature clearly embedded in the model.
\begin{table}[!ht]
\centering
\caption{Out-of-sample forecast accuracy measures RMSE, MAPE for the fitted ULARMA, KARMA and $\beta$ARMA.}\label{t:f}\vskip 0.5\baselineskip
\renewcommand{\arraystretch}{1.3}
\renewcommand{\tabcolsep}{2pt}
\footnotesize
\begin{tabular}{cccccccccccccc}
 \multicolumn{1}{c}{}& \multicolumn{12}{c}{Forecast Horizon ($h$)} \\
\cline{2-13}
\multirow{2}{*}{Model} & 1 & 2 & 3 & 4 & 5 & 6 & 7 & 8 & 9 & 10 & 11 & 12\\
\cline{2-13}
& \multicolumn{12}{c}{RMSE} \\
\hline
ULARMA  & 0.0015 & 0.0042 & 0.0042 & 0.0037 & 0.0052 & 0.0067 & 0.0070 & 0.0078 & 0.0083 & 0.0079 & 0.0077 & 0.0079\\ 
KARMA  & 0.0018 & 0.0014 & 0.0137 & 0.0156 & 0.0167 & 0.0165 & 0.0155 & 0.0148 & 0.0142 & 0.0135 & 0.0129 & 0.0124\\ 
$\beta$ARMA  & 0.0012 & 0.0015 & 0.0113 & 0.0126 & 0.0129 & 0.0139 & 0.0139 & 0.0142 & 0.0142 & 0.0135 & 0.0129 & 0.0124\\ 
\cline{2-13}
& \multicolumn{12}{c}{MAPE} \\
\hline
ULARMA  & 0.0242 & 0.0516 & 0.0575 & 0.0459 & 0.0698 & 0.0973 & 0.1078 & 0.1265 & 0.1391 & 0.1315 & 0.1273 & 0.1318\\ 
KARMA  & 0.0284 & 0.0207 & 0.1400 & 0.1827 & 0.2206 & 0.2362 & 0.2239 & 0.2194 & 0.2140 & 0.1947 & 0.1790 & 0.1718\\ 
$\beta$ARMA  & 0.0196 & 0.0222 & 0.1182 & 0.1491 & 0.1719 & 0.2038 & 0.2153 & 0.2317 & 0.2391 & 0.2203 & 0.2004 & 0.1920\\ 
\hline
\end{tabular}
\end{table}

Although stationarity conditions for GARMA and GARMA-like models based on continuous distributions are not known, it is customary to check for the presence of unitary roots in the model's AR part. The smallest root of the AR characteristic polynomial in absolute value (SRCP) is 1.13 for the ULARMA model, suggesting the absence of unit roots. For the $\beta$ARMA model, the SRCP is 1.054, slightly above the common 1.05 threshold often used as a rule of thumb, indicating a near-unit root.  In the case of the KARMA model, the SRCP values were 1.006 and 1.015, which are indicative of the presence of unit roots. This highlights the versatility of the ULARMA,  which was able to model the time series more parsimoniously than competitors, while avoiding the presence of unit roots.
\subsection{Discussion}
From the comparison of the proposed ULARMA model with the benchmark models KARMA and $\beta$ARMA for modeling the proportion of net electricity generated by conventional hydroelectric power in the United States, as presented in the previous section, several conclusions can be drawn. The empirical results emphasize the value of having a diverse repertoire of models available, and illustrate why it is essential to consider multiple alternatives when aiming for accurate forecasting.

While ULARMA relies on a uniparametric distribution, KARMA and $\beta$ARMA are based on two-parameter families, which could, in principle, offer greater flexibility. Nonetheless, in the presented application, ULARMA achieved a considerably more parsimonious fit relative to the benchmark models.

Regarding in-sample performance metrics, both KARMA and $\beta$ARMA consistently outperformed ULARMA. However, when the focus shifts to the most facet of time series analysis, out-of-sample forecasting, ULARMA yielded substantially better results, arguably due to its parsimonious structure. These findings highlight the adaptability of the Unit-Lindley distribution and, by extension, the practical effectiveness of the ULARMA model in real-world forecasting contexts.

\bibliographystyle{apalike}
\bibliography{bib}
\section*{Appendix A}
\begin{lema}\label{lema}
Let $Y\sim \mathrm{UL}(\mu)$, $\mu>0$. Then 
\begin{equation*}
\E\bigg(\frac{Y}{1-Y}\bigg)=\frac{\mu^2+\mu}{1-\mu}.
\end{equation*}
\end{lema}
\noindent \textbf{Proof:} Let $Y\sim \mathrm{UL}(\mu)$ and set $X:=Y/(1-Y)$. It is easy to see that $X$ has a probability density function given by
\begin{equation*}
    f_X(x;p,s)=\frac{1}{\mu}(1-\mu)^2(x+1)\exp\bigg\{\frac{x(\mu-1)}{\mu}\bigg\}I(x>0).
\end{equation*}
For brevity's sake, put $k:=(\mu-1)/\mu<0$. Then
\begin{align}\label{ex}
    \E(X)&=\frac{\mu}{k^2}\int_{0}^\infty t(t+1)e^{tk}dt=\frac{\mu}{k^2}\bigg[\int_0^\infty t^2e^{tk}dt+\int_0^\infty te^{tk}dt\bigg]\coloneqq\frac{\mu}{k^2}\bigg[A_1+A_2\bigg].
\end{align}
Using integration by parts, it follows that, for all $k<0$,
\begin{equation}\label{A2}
    A_2=\frac{te^{tk}}{k}\bigg |_{t=0}^\infty-\frac{1}{k}\int_0^\infty e^{tk}dt=\frac{1}{k^2},
\end{equation}
 which in turn implies that,
\begin{equation}\label{A1}
    A_1=\frac{t^2e^{tk}}{k}\bigg |_{t=0}^\infty-\frac{2}{k}\int_0^\infty te^{tk}dt=-\frac{2 A_2}{k}=-\frac{2}{k^3}.
\end{equation}
The result now follows by substituting \eqref{A2} and \eqref{A1} in \eqref{ex}.
\end{document}